\documentclass[11pt,reqno]{amsart}
\usepackage[margin=0.9in]{geometry}
\usepackage[utf8]{inputenc}
\usepackage[T1]{fontenc}
\usepackage{amssymb}
\usepackage{graphicx} 
\usepackage{mathrsfs}
\usepackage{enumerate}
\usepackage{xspace}
\usepackage{color}
\usepackage{enumerate}
\usepackage{enumitem}
\usepackage{amsthm}
\usepackage{dsfont}
\usepackage{bbm}
\usepackage{commath}
\usepackage[colorlinks=true,linkcolor=blue]{hyperref}
\usepackage[numbers]{natbib}
\usepackage[backgroundcolor=white,bordercolor=red]{todonotes}
\DeclareMathAlphabet{\mathpzc}{OT1}{pzc}{m}{it}
\usepackage[nameinlink]{cleveref}
\numberwithin{equation}{section}
\begin{document}

\expandafter\let\expandafter\oldproof\csname\string\proof\endcsname
\let\oldendproof\endproof
\renewenvironment{proof}[1][\proofname]{%
	\oldproof[\normalfont\scshape\hspace{1em}#1]%
}{\oldendproof}

\newtheoremstyle{mystyle_thm}
{6pt}
{6pt}
{\itshape}
{1em}
{\scshape}
{.}
{.5em}
{}%

\newtheoremstyle{mystyle_def}
{6pt}
{6pt}
{}
{1em}
{\scshape}
{.}
{.5em}
{}%

\theoremstyle{mystyle_thm}

\newtheorem{theorem}{Theorem}[section]
\newtheorem{lemma}[theorem]{Lemma}
\newtheorem{proposition}[theorem]{Proposition}
\newtheorem{corollary}[theorem]{Corollary}
\newtheorem{definition}[theorem]{Definition}
\newtheorem{Ass}[theorem]{Assumption}
\newtheorem{condition}[theorem]{Condition}
\newtheorem*{lemma*}{Lemma}

\theoremstyle{mystyle_def}

\newtheorem{example}[theorem]{Example}
\newtheorem{remark}[theorem]{Remark}
\newtheorem{SA}[theorem]{Standing Assumption}
\newtheorem{discussion}[theorem]{Discussion}
\newtheorem{remarks}[theorem]{Remark}
\newtheorem*{notation}{Remark on Notation}
\newtheorem{application}[theorem]{Application}

\newcommand{\of}{[\hspace{-0.06cm}[}
\newcommand{\gs}{]\hspace{-0.06cm}]}

\newcommand{\llambda}{{\mathchoice
		{\lambda\mkern-4.5mu{\raisebox{.4ex}{\scriptsize$\setminus$}}}
		{\lambda\mkern-4.83mu{\raisebox{.4ex}{\scriptsize$\setminus$}}}
		{\lambda\mkern-4.5mu{\raisebox{.2ex}{\footnotesize$\scriptscriptstyle\setminus$}}}
		{\lambda\mkern-5.0mu{\raisebox{.2ex}{\tiny$\scriptscriptstyle\setminus$}}}}}

\newcommand{\1}{\mathds{1}}

\newcommand{\F}{\mathbf{F}}
\newcommand{\G}{\mathbb{G}}

\newcommand{\B}{\mathbf{B}}

\newcommand{\M}{\mathcal{M}}

\newcommand{\la}{\langle}
\newcommand{\ra}{\rangle}

\newcommand{\lle}{\langle\hspace{-0.085cm}\langle}
\newcommand{\rre}{\rangle\hspace{-0.085cm}\rangle}
\newcommand{\blle}{\Big\langle\hspace{-0.155cm}\Big\langle}
\newcommand{\brre}{\Big\rangle\hspace{-0.155cm}\Big\rangle}

\newcommand{\X}{\mathsf{X}}

\newcommand{\lv}{I^-}
\newcommand{\uv}{I^+}

\newcommand{\tr}{\operatorname{tr}}
\newcommand{\N}{{\mathbb{N}}}
\newcommand{\cadlag}{c\`adl\`ag }
\newcommand{\on}{\operatorname}
\newcommand{\oP}{\overline{P}}
\newcommand{\oO}{\mathcal{O}}
\newcommand{\D}{\mathsf{D}} 
\newcommand{\bx}{\mathsf{x}}
\newcommand{\bb}{\hat{b}}
\newcommand{\bs}{\hat{\sigma}}
\newcommand{\bv}{\hat{v}}
\renewcommand{\v}{\mathfrak{m}}
\newcommand{\ob}{\bar{b}}
\newcommand{\oa}{\bar{a}}
\newcommand{\os}{\widehat{\sigma}}
\renewcommand{\j}{\varkappa}
\newcommand{\scl}{\ell}
\newcommand{\Y}{\mathscr{Y}}
\newcommand{\Z}{\mathscr{Z}}
\newcommand{\T}{\mathcal{T}}
\newcommand{\con}{\mathsf{c}}
\newcommand{\nk}{\hspace{-0.25cm}{{\phantom A}_k^n}}
\newcommand{\nl}{\hspace{-0.25cm}{{\phantom A}_1^n}}
\newcommand{\nm}{\hspace{-0.25cm}{{\phantom A}_2^n}}
\newcommand{\n}{\hspace{-0.35cm}{\phantom {Y_s}}^n}
\newcommand{\nme}{\hspace{-0.35cm}{\phantom {Y_s}}^{n - 1}}
\renewcommand{\o}{\hspace{-0.35cm}\phantom {Y_s}^0}
\newcommand{\e}{\hspace{-0.4cm}\phantom {U_s}^1}
\newcommand{\z}{\hspace{-0.4cm}\phantom {U_s}^2}
\newcommand{\iii}{|\hspace{-0.05cm}|\hspace{-0.05cm}|}
\newcommand{\co}{\overline{\on{co}}}
\newcommand{\ovb}{\overline{b}}
\newcommand{\ova}{\overline{a}}
\newcommand{\s}{\mathfrak{s}}
\newcommand{\opsi}{\overline{\Psi}}
\newcommand{\ol}{\mathcal{L}}
\newcommand{\cW}{\mathscr{W}}
\newcommand{\cU}{\mathcal{U}}
\newcommand{\oD}{\overline{D}}
\newcommand{\ua}{\underline{a}}
\newcommand{\ou}{\overline{b}}
\newcommand{\uu}{\underline{b}}
\renewcommand{\d}{\text{\rm d}}
\newcommand{\J}{J}
\renewcommand{\epsilon}{\varepsilon}

\newcommand{\fPs}{\fP_{\textup{sem}}}
\newcommand{\fPas}{\mathfrak{S}_{\textup{ac}}}
\newcommand{\rrarrow}{\twoheadrightarrow}
\newcommand{\cA}{\mathcal{A}}
\newcommand{\ocA}{\mathcal{U}}
\newcommand{\cR}{\mathcal{K}}
\newcommand{\cK}{\mathcal{K}}
\newcommand{\cQ}{\mathcal{Q}}
\newcommand{\cF}{\mathcal{F}}
\newcommand{\cE}{\mathcal{E}}
\newcommand{\cC}{\mathcal{C}}
\newcommand{\cD}{\mathcal{D}}
\newcommand{\bC}{\mathbb{C}}
\newcommand{\cH}{\mathcal{H}}
\newcommand{\bth}{\overset{\leftarrow}\theta}
\renewcommand{\th}{\theta}
\newcommand{\cG}{\mathcal{G}}
\newcommand{\fPasn}{\mathfrak{S}^{\textup{ac}, n}_{\textup{sem}}}
\newcommand{\CLM}{\mathfrak{M}^\textup{ac}_\textup{loc}}
\newcommand{\Sd}{\mathcal{S}^\textup{sp}_{\textup{d}}}
\newcommand{\Sc}{\mathcal{S}}
\newcommand{\Sac}{\mathcal{S}_\textup{ac}}
\newcommand{\A}{\mathbb{A}}
\newcommand{\Td}{\mathsf{T}^\textup{d}}
\renewcommand{\t}{\mathfrak{t}}

\newcommand{\bR}{\mathbb{R}}
\newcommand{\nnabla}{\nabla}
\newcommand{\f}{\mathfrak{f}}
\newcommand{\g}{\mathfrak{g}}
\newcommand{\oconv}{\overline{\on{co}}\hspace{0.075cm}}
\renewcommand{\a}{\mathfrak{a}}
\renewcommand{\b}{\mathfrak{b}}
\renewcommand{\d}{\mathsf{d}}
\newcommand{\bS}{\mathbb{S}^d_+}
\newcommand{\p}{\mathsf{p}}
\newcommand{\dr}{r} 
\newcommand{\m}{\mathbb{M}}
\newcommand{\Q}{Q}
\newcommand{\usc}{\textit{USC}}
\newcommand{\lsc}{\textit{LSC}}
\newcommand{\q}{\mathfrak{q}}
\renewcommand{\X}{\mathbb{X}}
\newcommand{\W}{\mathscr{W}}
\newcommand{\fP}{\mathcal{P}}
\newcommand{\w}{\mathsf{w}}
\newcommand{\oM}{\mathsf{M}}
\newcommand{\oZ}{\mathsf{Z}}
\newcommand{\oK}{\mathsf{K}}
\renewcommand{\Re}{\operatorname{Re}}
\newcommand{\cCk}{\mathsf{c}_k}
\newcommand{\C}{\mathsf{C}}
\newcommand{\cP}{\mathcal{P}}
\newcommand{\oPi}{\overline{\Pi}}
\newcommand{\cI}{\mathcal{I}}
\renewcommand{\P}{\mathbb{P}}
\newcommand{\E}{E}
\newcommand{\ug}{\underline{g}}

\newcommand{\UC}{\textit{UC}}
\renewcommand{\N}{\mathcal{C}}
\renewcommand{\G}{\mathbb{S}}
\renewcommand{\L}{{\mathcal{L}_\eta^2}}

\renewcommand{\emptyset}{\varnothing}

\allowdisplaybreaks

\makeatletter

 \title[Extremal solutions and stochastic control]{Stochastic Control Representations for Extremal Solutions of HJB Dirichlet Problems with Measurable Coefficients} 
\author[D. Criens]{David Criens}
\address{University of Freiburg, Ernst-Zermelo-Str. 1, 79104 Freiburg, Germany.}
\email{david.criens@stochastik.uni-freiburg.de}

\author[H. Park]{Hyungbin Park}
\address{Seoul National University, 1, Gwanak-ro, Gwanak-gu, Seoul, Republic of Korea.}
\email{hyungbin@snu.ac.kr, hyungbin2015@gmail.com}

\author[J. Park]{Jongjin Park}
\address{Seoul National University, 1, Gwanak-ro, Gwanak-gu, Seoul, Republic of Korea.}
\email{pjj4230@snu.ac.kr}



\thanks{David Criens gratefully acknowledges support from the Freiburg Center for Data Analysis, Modeling and AI. 
	Hyungbin Park was supported by National Research Foundation of Korea grants funded by the Ministry of Science and ICT (2022R1A5A6000840, RS-2026-25488333). }
\date{\today}

\begin{abstract}
	We study stochastic control problems for uniformly elliptic controlled diffusions on bounded domains whose coefficients are merely measurable in the state variable and whose drift is dominated by an \(L_d\)-function. The associated Hamilton--Jacobi--Bellman (HJB) Dirichlet problems may have more than one \(L_p\)-viscosity solution, so that the equation alone does not determine the value function. We show that the control problem selects the minimal solution: for all \(p\) in a range \((d_0,d)\) determined by the spatial dimension, the ellipticity constant and the drift bound, the value function is the minimal \(L_p\)-viscosity supersolution of the HJB Dirichlet problem. In particular, every continuous \(L_p\)-viscosity supersolution bounds the optimal cost from above, which yields a one-sided verification theorem without regularity assumptions. We further prove the existence of optimal relaxed control rules for running costs dominated by an \(L_p\)-function, and the lower semicontinuity of the value function under perturbations of the coefficients. An example shows that upper semicontinuity may fail even when the coefficients are Lipschitz continuous in the control variable and the approximating coefficients are smooth, so that regularization may overestimate the optimal cost. The proof approximates arbitrary supersolutions from below by Sobolev supersolutions of HJB equations whose coefficients are modified only on sets of vanishing measure, and passes to the limit by means of Krylov estimates
	and martingale problem techniques.

   
\end{abstract}

\maketitle

\section{Introduction} 
A cornerstone of stochastic optimal control is the characterization of the value function through the associated Hamilton--Jacobi--Bellman (HJB) equation. Under regularity assumptions on the coefficients, the value function is the unique viscosity solution of the HJB equation, which is classically shown either through the verification approach or via a comparison principle, see, e.g., \cite{CIL92, Lions85_2} and \cite{FS_06, pham} for textbook treatments. Both routes may break down when the coefficients are merely measurable in the state variable. A comparison principle is not available in general, the HJB equation need not have a Sobolev solution, and, most importantly, uniqueness may fail: the equation may have several, or even infinitely many, solutions (see \cite{zbMATH01140123, zbMATH01308655} and Appendix~\ref{ex: non unique}). Uniqueness can be restored under additional regularity assumptions on the diffusion coefficient. Under a uniform continuity assumption, \cite{deFeo_26} showed, for infinite time-horizon problems on \(\bR^d\) and \(p \geq d\), that \(L_p\)-viscosity solutions of the HJB equation are strong Sobolev solutions, and identified the value function as the unique \(L_p\)-viscosity solution by a verification argument. Without any continuity of the diffusion coefficient in the state variable, \cite{C_25_SIMA} showed, using mainly probabilistic arguments, that value functions are still \(L_p\)-viscosity solutions of their HJB equations. Since the equation alone no longer determines the value function, the following question arises:
\begin{enumerate}[label=(\roman*), ref=\roman*, series=questions]
	\item \label{qu: i} Which solution of the HJB equation does the control problem select?
\end{enumerate}

\noindent
Closely related are the further questions:

\begin{enumerate}[resume*=questions]
	\item \label{qu: ii} How does the lack of uniqueness affect the stability of the
	value function under perturbations of the coefficients, for instance under
	regularization?
	\item \label{qu: iii} Do optimal control rules exist when, in addition, the running
	cost is merely measurable in the state variable and possibly unbounded?
\end{enumerate}

Such questions arise naturally in applications. Coefficients that are merely measurable in the state variable describe dynamics that switch between regimes when the state crosses a threshold, as in advertising models considered in \cite{deFeo_26}. In stochastic differential games and mean-field problems, the coefficients faced by a single agent contain the feedback strategies of the other agents, which are in general only measurable, see \cite{zbMATH06441180, LakSPA15}. Finally, in computations the coefficients are routinely regularized or discretized. Without uniqueness for
the HJB equation, it is a priori unclear which solution, if any, such approximations recover. Question~\eqref{qu: ii} is motivated by this last point.

For a class of uniformly elliptic HJB Dirichlet problems with coefficients merely measurable in their space variable and possibly unbounded, but \(L_d\)-dominated, drift coefficient, we answer these questions as follows: \eqref{qu: i} the control problem selects the \emph{minimal} solution, more precisely, the value function is the minimal \(L_p\)-viscosity supersolution of the HJB Dirichlet problem (Theorem~\ref{theo: main1}), \eqref{qu: ii} the value function is lower semicontinuous, but in general not continuous, with respect to perturbations of the coefficients (Proposition~\ref{coro: lsc} and Remark~\ref{rem: lsc}), and \eqref{qu: iii} optimal relaxed control rules exist for running costs dominated by an \(L_p\)-function (Proposition~\ref{prop: optimal}). 

In the remainder of this introduction, we explain our setting and main results in a more precise way, comment on the proof strategy for our main result, and discuss related literature. 

\subsection{Setting and the non-uniqueness phenomenon}
Formally, we study the value function
\begin{align} \label{eq: v intro}
	v(x) := \inf_{\lambda, \, X^\lambda} E \Big[ f \big(X^\lambda_{\tau_D (X^\lambda)} \big) + \int_0^{\tau_D (X^\lambda)} g (\lambda_s, X^\lambda_s) \, ds \mid X^\lambda_0 = x\Big], \quad x \in \overline{D} = D \cup \partial D, 
\end{align} 
associated with the controlled process
\[
d X^\lambda_t = b (\lambda_t, X^{\lambda}_t) \, dt + \sqrt{a} (\lambda_t, X^\lambda_t) \, d W_t, \quad \tau_D (X^\lambda) := \inf \{t \geq 0 \colon X^\lambda_t \not \in D \}, 
\] 
where the infimum is taken over admissible relaxed control rules, i.e., joint laws of the control and the controlled process, see \cite{EKNJ88,nicole1987compactification,ElKa15, HausLep90}. The precise formulation is \eqref{eq: value function D} below. The diffusion coefficient \(a\) is
bounded and uniformly elliptic, the drift \(b\) is dominated by an
\(L_d\)-function, and both are continuous in the control variable \(\lambda\)
but merely Borel measurable in the state variable \(x\). The domain \(D\) is
bounded and satisfies the uniform exterior ball condition, the boundary
datum \(f\) is continuous, and the running cost \(g\) is continuous in \(\lambda\)
with \(\sup_{\lambda} |g (\lambda, \cdot \,)| \in L_p (D)\) for some \(p \in (d_0, d)\),
where \(d_0 \in (d/2, d)\) depends only on the dimension, the ellipticity
constant and the \(L_d\)-norm of the drift bound. The exponent \(d_0\) is provided by Krylov's theory for linear diffusions and elliptic equations with \(L_d\)-dominated drift \cite{krylov_CPDE_20, krylov_AOP_21,krylov_21}, and the restriction \(p < d\) comes from the existence theory in \cite{krylov_CPDE_20} and the \(L_p\)-viscosity theory for equations with \(L_d\)-drift from \cite{KS_22}, on which some of our proofs rely. For the critical case \(p = d\), the ABP maximum principle, a cornerstone of this theory, has recently been established in \cite{arXiv:2610.03918} under a slightly stronger condition on the drift. The relaxed control formulation provides the compactness that we use both for the existence of optimal controls and for all limiting arguments. For connections to other control frameworks, see~\cite{ElKa15}.

The associated HJB Dirichlet problem is
 \begin{align} \label{eq: DP intro}
 \begin{cases} \displaystyle - \inf_{\lambda \in \Lambda} \Big\{ \langle b (\lambda, x), \nabla u (x) \rangle + \tfrac{1}{2} \on{tr} \big[ a (\lambda, x) \nabla^2u(x) \big]+ g (\lambda, x)  \Big\} = 0, & x \in D, \\ u (x) = f (x), & x \in \partial D, \end{cases}
 \end{align} 
which belongs to the larger class of nonlinear PDEs that may not have classical solutions. 
We work with \(L_p\)-viscosity solutions, which have been introduced in the seminal work \cite{CCKS_96} and whose theory for HJB equations and more general nonlinear PDEs with measurable coefficients has since been developed extensively, see, e.g., \cite{JKS_02, Jen_Sw_05, KS_09, KS_22, KST_19, krylov_18} and the references therein.

Two features of this setting are important. First, the Dirichlet problem \eqref{eq: DP intro} need not have a Sobolev solution in \(W^2_{p, \on{loc}} (D) \cap C (\oD)\), since, by \cite[Proposition~4.5]{KS_22}, the existence of such a solution would imply uniqueness among \(L_p\)-viscosity solutions. Consequently, verification arguments
cannot be applied to \eqref{eq: DP intro} directly. Second, uniqueness does fail in general. Nadirashvili \cite{zbMATH01140123} constructed a bounded uniformly elliptic coefficient \(a \colon \bR^d \to \mathbb{S}^d_+\) such that, for \(d \geq 3, D = \{ x \in \bR^d \colon \|x \| < 1\}\) and suitable \(f \in C^2 (\partial D)\), the linear Dirichlet problem 
 \begin{align} \label{eq: linear example}
 \begin{cases}\on{tr} \big[ a (x) \nabla^2 u (x) \big] = 0, & \|x\| < 1, \\ 
 	u (x) = f (x), & \|x \|  = 1, 
 \end{cases} 
 \end{align} 
 has at least two \(L_p\)-viscosity solutions (here, we use the fact that ``good'' solutions as studied in \cite{zbMATH01140123} are \(L_p\)-viscosity solutions, see \cite{CCKS_96, JKS_02}), see also \cite{zbMATH01308655} for a related construction. In Appendix~\ref{ex: non unique}, we add a small controlled drift to this example and obtain a semilinear HJB Dirichlet problem of the form \eqref{eq: DP intro}, nonlinear in the gradient, with at least two, and in fact infinitely many, \(L_p\)-viscosity solutions.

  The solution set nevertheless has some structure. For Hamiltonians with bounded coefficients,~\cite{Jen_Sw_05} proved that minimal and maximal \(L_p\)-viscosity solutions exist and, for the linear problem \eqref{eq: linear example}, \cite[Discussion~2.10 (b)]{C_25_SIMA} represented both stochastically, as infimum and supremum over all solutions of the associated martingale problem. For bounded and lower semicontinuous running costs, \cite{C_25_SIMA} also showed that the value function \(v\) is an \(L_p\)-viscosity solution of \eqref{eq: DP intro} for every \(p\geq d_0\).
 
\subsection{Main results}

Our main result, Theorem~\ref{theo: main1}, states that, for every \(p \in (d_0, d)\),
\[
v = \min \big\{ u\in C(\overline{D}) \colon u \text{ is an } L_p\text{-viscosity supersolution of
	\eqref{eq: DP intro}} \big\}.
\]
In words, \(v\) is an \(L_p\)-viscosity supersolution of the Dirichlet problem \eqref{eq: DP intro}, now
also for running costs that are only dominated by an \(L_p\)-function, and \(v \leq u\) for every \(L_p\)-viscosity supersolution \(u\). We emphasize three aspects of this result.

Among all solutions of \eqref{eq: DP intro}, the control problem selects the minimal one, and \(v\)
is given by a Perron-type formula \cite{zbMATH04142546, zbMATH05269301} as the pointwise minimum of all supersolutions. In particular, Theorem~\ref{theo: main1} provides the existence of a minimal \(L_p\)-viscosity solution for HJB equations with unbounded, \(L_d\)-dominated drift, extending the existence result of \cite{Jen_Sw_05} for bounded coefficients within the class of HJB equations considered here. The Perron method from \cite{zbMATH05269301} assumes locally bounded coefficients, which exclude an \(L_d\)-dominated drift.
From the control point of view, Theorem~\ref{theo: main1} can be viewed as part of a one-sided verification theorem for \(L_p\)-viscosity solutions that requires no regularity. Every \(L_p\)-viscosity supersolution of \eqref{eq: DP intro} is an
upper bound for the optimal cost, and \(v\) is the best such bound. This complements verification theorems for viscosity solutions \cite{zbMATH02205799, zbMATH05822711, zbMATH01000474} and the stochastic Perron method~\cite{BS_13}.
Theorem~\ref{theo: main1} also extends the stochastic representation of the minimal solution from the linear setting of \cite{C_25_SIMA} to nonlinear HJB equations. By Proposition~\ref{prop: representation}, it covers Dirichlet problems for concave Hamiltonians that differ by at most an \(L_p\)-function from a positively homogeneous, concave, uniformly elliptic Hamiltonian with \(L_d\)-drift bound. Applied to the costs \(-f\) and \(-g\), it represents the maximal \(L_p\)-viscosity subsolution of the maximization problem
 \begin{align} \label{eq: sup-DP}
 	\begin{cases} \displaystyle - \sup_{\lambda \in \Lambda} \Big\{ \langle b (\lambda, x), \nabla u (x) \rangle + \tfrac{1}{2} \on{tr} \big[ a (\lambda, x) \nabla^2u(x) \big] + g (\lambda, x)\Big\} = 0, & x \in D, \\ u (x) = f (x), & x \in \partial D, \end{cases}
 \end{align}
by the supremum over admissible relaxed control rules, see Corollary~\ref{coro: main1}. This supremum need not solve the minimization problem~\eqref{eq: DP intro}.

Classical results identify the value function of a minimization problem as the maximal subsolution of its HJB equation in a weak sense \cite{Lions85_1, LM80, zbMATH03752736}, or as the supremum of its smooth subsolutions \cite{zbMATH04135852}. Under non-uniqueness, a characterization as in \cite{Lions85_1, LM80, zbMATH03752736} is not available. Since \(v\) is the minimal solution, it can be the maximal subsolution only if \eqref{eq: DP intro} has a unique solution, see Remark~\ref{rem: asymmetry}~(a). The reason for this asymmetry is structural. While the proof of \(v \leq u\) for a supersolution \(u\) is related to tightness and closedness properties, the proof of \(w \leq v\) for a subsolution \(w\) may fail without a uniqueness property of a controlled martingale problem, see Remark~\ref{rem: asymmetry}~(b) for a discussion. A related asymmetry was observed in \cite{nicole1987compactification} for the semicontinuity of value functions. To the best of our knowledge, it remains an interesting open problem whether \eqref{eq: DP intro} admits a maximal \(L_p\)-viscosity solution in the present setting, and, if so, whether it has an intrinsic stochastic representation.

The same asymmetry is reflected in the stability question. For coefficient perturbations satisfying the hypotheses of Proposition~\ref{coro: lsc}, the associated value functions \(v^{(n)}\) satisfy
\begin{align*}
	v(x) \leq \liminf_{n \to \infty} v^{(n)} (x),\quad x \in \oD.
\end{align*}
In particular, if the coefficients are regularized, for instance mollified in the state variable, the regularized value functions can, in the limit, only overestimate the optimal cost, see Remark~\ref{rem: lsc}~(a). The example in Appendix~\ref{ex: non unique} shows that the overestimation can be strict, even when the
coefficients are Lipschitz continuous in the control variable and the approximating coefficients are smooth, see Remark~\ref{rem: lsc}~(b). Related one-sided stability results under drift and cost mollification for parabolic finite time-horizon problems with additive noise were obtained in \cite[Section~3.1]{arXiv:2509.01054}, see Remark~\ref{rem: lsc}~(b) for a comparison. Finally,
Proposition~\ref{prop: optimal} provides optimal relaxed control rules for running costs
dominated by an \(L_p\)-function, complementing \cite{C_25_SIMA}, where optimal
strong Markov control rules are obtained for bounded lower semicontinuous
running costs.

\subsection{Strategy of the proof}

Since an arbitrary \(L_p\)-viscosity supersolution \(u\) need not be regular
enough for a verification argument, and no comparison principle is
available, we approximate \(u\) by functions to which verification applies.
In the linear case without first-order terms and with bounded right-hand
side, \cite{J_96} showed that \(L_p\)-viscosity solutions are limits of
strong solutions of linear equations with approximating coefficients, see also \cite{JKS_02}.
To treat the linear equation \eqref{eq: linear example}, \cite{C_25_SIMA} combined this with probabilistic limit arguments. The same general pattern underlies the control representations of
extremal viscosity solutions for equations with (semi)continuous
coefficients in \cite{zbMATH05032612, PCP_26, zbMATH01493589}, with approximations of a different
nature.

In our nonlinear setting, the approximation must
preserve the control structure, the ellipticity constant and the
\(L_d\)-drift bound, and limits of the associated control rules must be
admissible for the original problem.

Our main technical contribution, Lemma~\ref{lem: main approx}, provides such an approximation for drifts that are only \(L_d\)-dominated. For an arbitrary \(L_p\)-viscosity supersolution \(u\) and an exhaustion \((D_m)_{m=1}^\infty\) of \(D\), it yields Sobolev functions \(u_m \le u\) with \(\sup_{D_m}|u_m - u| \to 0\) that are supersolutions of modified HJB equations, whose coefficients satisfy the original ellipticity and
\(L_d\)-drift bounds and coincide with the original coefficients outside sets \(G_m\) whose Lebesgue measure tends to zero. The construction combines a penalized Pucci extremal equation with contact-set arguments, adapting ideas from \cite{zbMATH03751685, zbMATH05032612, JKS_02, Jen_Sw_05, J_96}, the solvability theory for nonlinear equations with \(L_d\)-drift \cite{krylov_CPDE_20}, and a Krylov--Safonov type estimate \cite{krylov_AOP_21, zbMATH03726824}, see Remark~\ref{rem: comparison} for more comments on the relation to these references. A measurable selection of minimizers then yields feedback controls whose costs are bounded by \(u_m\), which follows from a verification argument with Krylov's It\^o formula for Sobolev functions \cite{krylov_AOP_21}, see  Lemma~\ref{lem: stochastic lower bound}. Tightness, Krylov estimates and the martingale problem formulation of relaxed controls allow us to pass to the limit and to conclude \(v \leq u\), see Lemmas~\ref{lem: tight + closed} and~\ref{lem: continuity}.

For bounded drift coefficients, the approximation step could alternatively
be based on the proof of \cite[Theorem~3.1]{JKS_02}, together with a Sobolev existence result such as \cite[Theorem~1.2]{krylov_CPDE_20}. The \(L_d\)-dominated case
requires the new construction outlined above.

\subsection{Related literature}
Krylov estimates and It\^o's formula for Sobolev functions, the basic tools for controlled non-degenerate diffusions without smoothness, go back to \cite{krylov_80}.
More recently, stochastic control problems with measurable coefficients have been studied in \cite{C_25_SIMA, deFeo_26, arXiv:2509.01054, MPT_23}. The paper \cite{MPT_23} develops a maximum principle for controlled stochastic differential equations with additive noise and measurable drift, while \cite{deFeo_26,arXiv:2509.01054} use the verification approach to relate value functions to HJB equations. There, \cite[Proposition~3.3]{deFeo_26} assumes that the diffusion coefficient is uniformly
continuous on bounded spatial sets, uniformly in the control, and \cite[Theorem~2.1]{arXiv:2509.01054} has a constant non-degenerate diffusion coefficient, in which case the parabolic HJB equation has a unique Sobolev solution, which is identified as the corresponding value function. In the present paper, the diffusion coefficient
is merely measurable, so that, as explained above, verification cannot be
applied to \eqref{eq: DP intro} itself. For relaxed controls and associated compactness methods, see \cite{EKNJ88, nicole1987compactification, ElKa15, HausLep90}. Stochastic representations through dynamic programming and comparison are developed in \cite{CN22b, CN_25_JEEQ, CN_25_EJP, FS_06, Lions85_2, pham}, while \cite{BS_13, BT_11} provide approaches that avoid an a priori proof of the dynamic programming principle and measurable selection arguments, respectively. 

Besides the classical characterizations \cite{zbMATH04135852, Lions85_1, LM80, zbMATH03752736} discussed above, extremal characterizations of value functions without uniqueness have been obtained for infinite-horizon problems with Lipschitz coefficients
\cite{zbMATH05032612}, where continuous approximations of a supersolution serve as
obstacles in discounted problems to which comparison \cite{ishii} applies, for
first-order equations \cite{zbMATH01493589} and for first-order equations with
discontinuous running costs \cite{zbMATH02113830, zbMATH05131146}. Obstacle problems for Hamiltonians mollified in the
state variable also underlie the construction of maximal and minimal
\(L_p\)-viscosity solutions in \cite[Theorem~4.1]{Jen_Sw_05}, which, however, does
not yield a stochastic representation. Without convexity of the
Hamiltonian, \cite{zbMATH06261003} shows for parabolic equations with bounded
measurable coefficients that solutions of cutoff equations converge to the
maximal \(L_{d+1}\)-viscosity solution, see \cite{zbMATH06249180} for elliptic cutoff
equations and \cite{zbMATH06288066} for their identification with value functions of
stochastic differential games under Lipschitz conditions.

The present paper builds on \cite{C_25_SIMA}, which treats the linear problem \eqref{eq: linear example}, where both extremal solutions are represented, and shows the viscosity property of \(v\) for bounded lower semicontinuous running costs and continuity of \(v\)
with respect to the coefficients under an oscillation condition ensuring
uniqueness, see \cite[Proposition~2.9]{C_25_SIMA}. New here are the nonlinear minimal
representation, running costs dominated by an \(L_p\)-function, one-sided
stability without any condition ensuring uniqueness, and the example showing
that two-sided stability fails. In the companion paper \cite{PCP_26}, the value
function of parabolic terminal-payoff problems with killing is identified as
the greatest bounded upper semicontinuous \(C\)-viscosity subsolution, for
jointly upper semicontinuous sublinear Hamiltonians with locally bounded
coefficients under a Lyapunov condition \cite[Theorem~5.3]{PCP_26}. This may be
viewed as a parabolic analog of Corollary~\ref{coro: main1} under different assumptions.
Its proof approximates the Hamiltonian from above by continuous Hamiltonians, followed by a regularization and a non-degenerate perturbation. Such an approximation is not available here, since pointwise recovery of a Hamiltonian from continuous majorants would force upper semicontinuity in
the state variable. Lemma~\ref{lem: main approx} replaces it by an approximation of the supersolution itself.

\subsection{Organization} 

Section~\ref{sec: setting} introduces the relaxed control framework. Section~\ref{sec: main results} contains the main results: Theorem~\ref{theo: main1} with its consequences (Section~\ref{sec: main result}), stability
with respect to the coefficients (Section~\ref{sec: stability}), and the class of concave
Hamiltonians covered by our setting (Section~\ref{sec: concave}). Section~\ref{sec: proof} provides the proof of Theorem~\ref{theo: main1}. Appendix~\ref{ex: non unique} constructs a semilinear HJB Dirichlet problem with infinitely many \(L_p\)-viscosity solutions, and Appendix~\ref{app: pf prop} proves
Proposition~\ref{prop: representation}.

\section{The Control Setting} \label{sec: setting}

Let \(\Lambda\) be a non-empty compact Polish space and fix a spatial dimension \(d \geq 2\).
The set of all symmetric \(d \times d\) matrices with real entries is denoted by \(\mathbb{S}^d\) and the set of all positive semidefinite \(A \in \mathbb{S}^d\) is denoted by \(\bS\). The coefficients of the control problem under investigation are the two functions 
\begin{align*}
	a \colon\Lambda \times \bR^d \to \bS, \qquad
	b \colon \Lambda \times \bR^d \to \bR^d.
\end{align*}
Throughout this paper, we impose the following assumptions on the coefficients. 

\begin{SA} \label{SA: main1}
	\begin{enumerate}
		\item[(a)] \(a, b\) are Borel functions and, for every compact set \(K \subset \bR^d\), there exists a modulus of continuity\footnote{i.e., an increasing function with \(w_K (x) \to 0\) as \(x \searrow 0\)} \(w_K \colon \bR_+ \to \bR_+\) and a non-negative Borel function \(\mathfrak{c}_K \in L_d (K)\) such that 
		\[
		\| a (\lambda, x) - a (\lambda', x) \| + \| b (\lambda, x) - b (\lambda', x) \| \leq \mathfrak{c}_K (x) \, w_K \big(d_\Lambda (\lambda, \lambda') \big) 
		\] 
		for all \(\lambda, \lambda' \in \Lambda\) and \(x \in K\). Here, \(d_\Lambda\) denotes a metric on \(\Lambda\) that induces its topology. 
		\item[(b)] \(a\) is bounded and uniformly elliptic in the sense that there exists a \(\delta \in (0, 1)\) such that
		\[
		a (\lambda, x) \in \mathbb{S}_\delta^d := \Big\{ A \in \mathbb{S}^d_+ \colon \delta \, \|\xi\|^2 \leq \langle \xi, A \xi \rangle \leq \frac{\|\xi\|^2}{\delta}, \ \ \forall \, \xi \in \bR^d\Big\}
		\] 
		for all \((\lambda, x) \in \Lambda \times \bR^d\).
		\item[(c)] There exists a Borel function \(\b \in L_d (\bR^d)\) such that \(\|b (\lambda, x)\| \leq \b (x)\) for all \(\lambda \in \Lambda\) and \(x \in \bR^d\).
	\end{enumerate}
\end{SA}

\begin{example} \label{ex: coefficients}
	An example one may keep in mind is the following: 
	For \(\delta \in (0, 1)\), take \(\Lambda := [- 1, 1]\), \(\alpha \in (0, 1 - \delta)\) and
	\[
	a (\lambda, x) := (1 - \alpha \lambda) A (x), \quad b (\lambda, x) := \lambda h (x) e, 
	\] 
	where \(A \colon \bR^d \to \mathbb{S}^d_{\delta / (1 - \alpha)}\) is a Borel map, \(e \in \bR^d\) and \(h \in L_d (\bR^d)\). 
\end{example}

Let \(\Omega\) be the space of continuous functions \(\omega \colon \bR_+ \to \bR^d\) and endow it with the local uniform topology. The coordinate process on \(\Omega\) is denoted by \(X\), i.e., \(X_t (\omega) = \omega (t)\) for \(\omega \in \Omega\) and \(t \in \bR_+\). Further, \(\cF\) and \((\cF_t)_{t \geq 0}\) are the \(\sigma\)-field and filtration generated by \(X\), i.e., \(\cF := \sigma (X_t \colon t \geq 0)\) and \(\cF_t := \sigma (X_s \colon s \leq t)\).
Next, we introduce the space of control rules, which are formulated in the relaxed framework (see \cite{EKNJ88,nicole1987compactification, ElKa15}).
Let \(\m\) be the set of all Radon measures on \(\bR_+ \times \Lambda\) whose projections to \(\bR_+\) coincide with the Lebesgue measure. 
We endow \(\m\) with the vague topology, which turns it into a compact Polish space (\cite[Theorem 2.2]{EKNJ88}). The identity map on \(\m\) is denoted by \(M\), i.e.,
\[
M (d\lambda, ds) (m) := m (d \lambda, ds), \quad m \in \m.
\]
Define the \(\sigma\)-field \[\M := \sigma (M_t (\phi) \colon t \in \bR_+, \phi \in C_{c} (\bR_+ \times \Lambda)),\] where
\[
M_t (\phi) := \int_0^t \int_{\Lambda} \phi (s, \lambda) \, M(d\lambda, ds).
\]
At this point, we emphasize that \((M_t (\phi))_{t \geq 0}\) is a stochastic process on the space \(\m\).
In the following, the product space \((\Omega \times \m, \cF \otimes \M)\) will be of interest. Extending the notation above, we denote the coordinate process on this space by \((X, M)\). 
For \(\varphi \in C^2_c (\bR^d)\), define
\begin{align} \label{eq: L}
	L (\lambda, x, \nabla \varphi (x), \nabla^2 \varphi (x)) &:=  \langle b (\lambda, x), \nabla \varphi (x) \rangle + \tfrac{1}{2} \on{tr} \big[ a (\lambda, x) \nabla^2 \varphi (x) \big], \quad x \in \bR^d,
\end{align}
and 
\begin{align} \label{eq: test process}
	C_t (\varphi) := \varphi (X_t) - \varphi (X_0) - \int_{0}^{t} \int_{\Lambda}  L (\lambda, X_s, \nabla \varphi (X_s), \nabla^2 \varphi (X_s)) \, M (d \lambda, ds), \quad t \in \bR_+.
\end{align} 
A \emph{relaxed control rule} with initial value \(x\in\bR^d\) is a probability measure \(P\) on the product space \((\Omega \times \m, \mathcal{F} \otimes \M)\) with \(P (X_0 = x) = 1\) such that the processes \((C_t (\varphi))_{t \geq 0}\), with \(\varphi \in C^2_c (\bR^d)\), are \(P\)-martingales for the filtration
\[
\cG_s := \sigma (X_r, M_r (\phi) \colon r \leq s, \phi \in C_c (\bR_+ \times \Lambda)), \quad s \in \bR_+.
\]
\begin{remark}
	By the boundedness hypothesis on the coefficient \(a\), without loss of generality, we may also replace the class \(C^2_c (\bR^d)\) of test functions by \(C^2_b (\bR^d)\), see the proof of \cite[Proposition~5.4.11]{KaraShre} for details. 
\end{remark}
For \(x \in \bR^d\), we write 
\begin{align*}
\cR	(x)&:= \big\{ \text{all relaxed control rules with initial value } x \big \}.
\end{align*}

As the following example illustrates, the \(L_d\)-drift bound from Standing Assumption~\ref{SA: main1}~(c) is sharp within the scope of \(L_p\)-conditions that ensure that \(\cR (x) \not = \emptyset\).

\begin{example} \label{ex: existence}
	Adapting \cite[Example~1.1]{krylov_AOP_21} to our controlled setting, in case \(a \equiv \on{id}_d\) is the \(d\)-dimensional unit matrix and \(b (\lambda, x) \equiv b (x) := - x \| x\|^{- 2} (d / 2) \1_{\{ 0 < \|x\| < 1\}}\), we have \(\cR (0) = \emptyset\). It is easy to see that \(b \in L_{d - \varepsilon} (\bR^d)\) for all \(\varepsilon \in (0, d)\), but \(b \not \in L_d (\bR^d)\), which explains the above optimality claim. 
\end{example}

\section{A Stochastic Control Representation for Extremal Solutions}  \label{sec: main results}

\subsection{Main results} \label{sec: main result}

Let \(D \subset \bR^d\) be a non-empty bounded domain. As usual, the boundary of \(D\) is denoted by \(\partial D := \overline{D} \setminus D\), where \(\overline{D}\) denotes the closure in \(\bR^d\). Recall that \(D\) is said to satisfy the uniform exterior ball condition if there exists a radius \(\rho > 0\) such that, for every \(x \in \partial D\), there exists a center \(z \not \in \oD\) such that \(\overline{B}_{\rho} (z) \cap \oD= \{x\}\), where \(\overline{B}_{\rho} (z) := \{ y \in \bR^d \colon \|y - z\| \leq \rho\}\). 
In addition to the coefficients \(a\) and \(b\) from the previous section, we consider Borel functions
\begin{align*}
	f &\colon \partial D\to \bR, \quad
	g \colon \Lambda \times D \to \bR, 
\end{align*}
and impose the following:

\begin{SA} \label{SA: main2}  \(D\) satisfies the uniform exterior ball condition and \(f \in C (\partial D)\). 
Moreover, there exists a non-negative Borel function \(\mathfrak{g} \colon D \to \bR_+\) and a modulus of continuity \(w_g \colon \bR_+ \to \bR_+\) such that 
\[
\sup_{\lambda \in \Lambda} | g (\lambda, x)| \leq \g (x), \qquad | g (\lambda^*, x) - g (\lambda', x) | \leq \g (x) \, w_g \big( d_\Lambda (\lambda^*, \lambda') \big)
\] 
for all \(\lambda^*, \lambda' \in \Lambda\) and  \(x \in D\).
\end{SA}
We deduce from \cite[Lemmas~3.5, 3.6]{C_25_SIMA} that \(\cR (x)\) is a non-empty compact set (in the space of probability measures on \((\Omega \times \m, \cF \otimes \mathcal{M})\) with the weak topology) and that \(P\)-a.s. \(\tau_D < \infty\) for all \(P \in \cR(x)\), where 
\[
\tau_D := \inf \{t \geq 0 \colon X_t \not \in D \}.
\] 
Moreover, by 
\cite[Theorem~2.8]{krylov_AOP_21}, 
we have the following Krylov estimate: For every bounded domain \(G \subset \bR^d\), there exist two constants \[d_0 = d_0 (d, \delta, \| \b \|_{L_d (\bR^d)}) \in (d/2, d) \quad \text{and} \quad N = N (d, \delta, \| \b \|_{L_d (\bR^d)}, \on{diam} \, (G)) > 0\] such that, for all \(p \geq d_0\), every non-negative Borel function \(\ell \colon G \to [0, \infty]\) and \(x \in \bR^d\), 
\begin{align} \label{eq: krylov estimate} 
	\sup_{P \in \cK (x)} E^P \Big[ \int_0^{\tau_G} \ell (X_s) \, ds \Big] \leq N \| \ell \|_{L_p (G)}.
\end{align} 
In view of these facts, in case \(\sup_{\lambda \in \Lambda} |g (\lambda, \cdot \,)| \in L_p (D)\) for some \(p \geq d_0\), we may define
\begin{align} \label{eq: value function D}
	v (x) := \inf_{P \in \cR(x)} E^P & \Big[ f ( X_{\tau_D}) + \int_0^{\tau_D} \int_\Lambda g (\lambda, X_s) \, M (d \lambda, ds) \Big], \quad x \in \oD.
\end{align}
This function is called the {\em value function}. 
The value function is related to a Hamilton--Jacobi--Bellman (HJB) Dirichlet problem
\begin{align} \label{eq: PDE}
\begin{cases}	- H (x, \nabla u (x), \nabla^2 u (x)) = 0, & x \in D,  \\
	u (x) = f (x), & x \in \partial D, \end{cases}
\end{align} 
where, for \(x \in D, q \in \bR^d\) and \(A \in \mathbb{S}^d\), 
\begin{align} \label{eq: Hamiltonian}
	H (x,q, A) &:=  \inf_{\lambda \in \Lambda}  \Big\{ \langle b (\lambda, x), q\rangle + \tfrac{1}{2} \on{tr} \big[ a (\lambda, x) A \big] + g (\lambda, x)\Big\}.
\end{align}
As the coefficients of this equation are merely measurable, we use the concept of \(L_p\)-viscosity solutions as introduced in the seminal paper~\cite{CCKS_96}. For the reader's convenience, we recall the definition. Throughout, we always assume that \(p > d/2\) and identify each element of \(W^2_{p,\on{loc}} (D)\) with its unique continuous representative, whose existence follows from the Sobolev embedding theorem. For its weak derivatives of orders one and two, we fix Borel measurable representatives.

\begin{definition} \label{def: VS elliptic}
	A continuous function \(u \colon \oD \to \bR\) is called an {\em \(L_p\)-viscosity subsolution} (resp., {\em \(L_p\)-viscosity supersolution}) to \eqref{eq: PDE}
	if \(u \leq f\) on \(\partial D\) (resp. \(u \geq f\) on \(\partial D\)) and for every \(\varphi \in W^{2}_{p, \textup{loc}} (D)\) such that \(u - \varphi\) has a local maximum (resp., a local minimum) at some \(x_0 \in D\), we have 
	\begin{align*}
		&\lim\nolimits_{r \searrow 0} \on{ess \, inf}_{x \in B_r (x_0)} \, ( - H (x, \nabla \varphi (x), \nabla^2 \varphi (x)) ) \leq 0 \\
		(\text{resp. }\hspace{-0.05cm}&\lim\nolimits_{r \searrow 0} \on{ess \, sup}_{x \in B_r (x_0)} \,
		(- H (x, \nabla \varphi (x), \nabla^2 \varphi (x))) \geq 0), 
	\end{align*} 
	where \(B_r (x_0) := \{ y \in \bR^d \colon \|y - x_0\| < r \}\).
	The function \(u\) is an {\em \(L_p\)-viscosity solution} to \eqref{eq: PDE} if it is both an \(L_p\)-viscosity sub- and supersolution.
\end{definition}

If \(g\) were bounded and lower semicontinuous on \(\overline{D}\), \cite[Theorem~2.5]{C_25_SIMA} shows that there exists a constant \(d_0 = d_0 (d, \delta, \| \b\|_{L_d (\bR^d)}) \in (d/2, d)\) such that \(v\) is an \(L_p\)-viscosity solution to the equation \eqref{eq: PDE} for all \(p \geq d_0\). The next theorem extends this result to more general \(g\) and, more importantly, shows that \(v\) has a distinguished role within the class of \(L_{p}\)-viscosity supersolutions as being the minimal one. In particular, it proves the existence of a minimal \(L_p\)-viscosity solution within a setting with possibly unbounded but \(L_d\)-dominated drift. 

\begin{theorem} \label{theo: main1}
	We impose the Standing Assumptions~\ref{SA: main1} and \ref{SA: main2}. There exists a constant \(d_0 = d_0 (d, \delta, \| \b\|_{L_d (\bR^d)}) \in (d/2, d)\) such that, for every \(p \in (d_0, d)\) with \(\mathfrak{g} \in L_p (D)\), the value function \(v\) is an \(L_p\)-viscosity solution to the Dirichlet problem \eqref{eq: PDE}. Moreover, it is the minimal \(L_p\)-viscosity supersolution, i.e., if \(u\) is another \(L_p\)-viscosity supersolution, then \(v \leq u\) on \(\overline{D}\). 
\end{theorem}

The proof of Theorem~\ref{theo: main1} is given in Section~\ref{sec: proof} below.

\begin{corollary} \label{coro: main1}
	Assume the Standing Assumptions~\ref{SA: main1} and \ref{SA: main2} with \(\mathfrak{g} \in L_p (D)\) for some \(p \in (d_0, d)\) and \(d_0\) as in Theorem~\ref{theo: main1}. Then, the value function 
	\begin{align} \label{eq: u maxi}
	\overline{v} (x) := \sup_{P \in \cK (x)} E^P \Big[ f (X_{\tau_D}) + \int_0^{\tau_D} \int_\Lambda g (\lambda, X_s) \, M (d \lambda, ds)\Big], \quad x \in \overline{D},  
	\end{align} 
	is an \(L_p\)-viscosity solution to the Dirichlet problem
	\[
	\begin{cases} \displaystyle - \sup_{\lambda \in \Lambda} \Big\{ \langle b (\lambda, x), \nabla u (x) \rangle + \tfrac{1}{2} \on{tr} \big[ a (\lambda, x) \nabla^2u(x) \big] + g (\lambda, x)\Big\} =0, & x \in D, \\ u (x) = f (x), & x \in \partial D.
	\end{cases}
	\]
	Moreover, it is the maximal \(L_p\)-viscosity subsolution.  
\end{corollary}

\begin{proof}
	This follows from Theorem~\ref{theo: main1} applied with the negated functions \(- f\) and \(- g\). 
\end{proof}

As a sort of byproduct of our proof for Theorem~\ref{theo: main1}, we establish the existence of optimal control rules. This complements the somewhat stronger observation from \cite{C_25_SIMA} that optimal strong Markov control rules exist whenever \(g\) is bounded and lower semicontinuous. 

\begin{proposition} \label{prop: optimal}
	Assume the Standing Assumptions~\ref{SA: main1} and \ref{SA: main2} with \(\mathfrak{g} \in L_{d_0} (D)\) for \(d_0\) as in Theorem~\ref{theo: main1}. Then, there exists a family \((P_x)_{x \in D}\), with \(P_x \in \cK (x)\) for \(x \in D\), such that 
	\[
	v (x) = E^{P_x} \Big[ f (X_{\tau_D}) + \int_0^{\tau_D} \int_\Lambda g (\lambda, X_s) \, M (d \lambda, ds) \Big] 
	\] 
	for all \(x \in D\). 
\end{proposition}

\begin{proof} 
We recall from \cite[Lemma~3.6]{C_25_SIMA} that \(\cK (x)\) is compact in the weak topology. Moreover, continuity (in the weak topology) of 
\[
\cK (x) \ni P \mapsto E^P \Big[ f (X_{\tau_D}) + \int_0^{\tau_D} \int_\Lambda g (\lambda, X_s) \, M (d \lambda, ds) \Big]
\]
can be established analogously to Lemma~\ref{lem: continuity} below. These two facts entail the claim.
\end{proof}

		\begin{remark} \label{rem: asymmetry}
		Theorem~\ref{theo: main1} identifies \(v\) as the minimal \(L_p\)-viscosity solution of \eqref{eq: PDE}, but contains no statement on the maximal one.
		
		(a) If \(v\) were also the maximal \(L_p\)-viscosity subsolution of \eqref{eq: PDE}, then every \(L_p\)-viscosity solution \(u\) would satisfy \(v \leq u \leq v\), i.e., \(v\) would be the unique \(L_p\)-viscosity solution. As \eqref{eq: PDE} may have more than one \(L_p\)-viscosity solution (see \cite{zbMATH01140123} and Appendix~\ref{ex: non unique}), \(v\) is, in general, not the maximal \(L_p\)-viscosity subsolution. In the linear case, i.e., if \(\Lambda\) is a singleton, the maximal \(L_p\)-viscosity solution of \eqref{eq: PDE} is given by \eqref{eq: u maxi}, i.e., by the supremum over all solutions of the martingale problem. This follows directly from Theorem~\ref{theo: main1} and Corollary~\ref{coro: main1}. In the nonlinear case, however, \eqref{eq: u maxi} is the maximal \(L_p\)-viscosity subsolution of the \(\sup\)-problem \eqref{eq: sup-DP} and, as \(\sup_{\lambda \in \Lambda} \geq \inf_{\lambda \in \Lambda}\), only an \(L_p\)-viscosity supersolution of \eqref{eq: PDE}. In general, it does not solve \eqref{eq: PDE}. To the best of our knowledge, it remains an open problem whether, in our setting with \(L_d\)-dominated drift, \eqref{eq: PDE} always admits a maximal \(L_p\)-viscosity solution (for bounded coefficients, this follows from \cite{Jen_Sw_05}), and, if so, whether it admits an intrinsic stochastic representation.
		
		(b) On the level of the proof, the minimality of \(v\) only relies on tightness and closedness properties of control rules, while a maximality statement would require an approximation property of the controlled martingale problem. A related asymmetry was observed in \cite{nicole1987compactification} for the regularity of value functions: one semicontinuity property follows from tightness and closedness, whereas the other one is related to uniqueness for the controlled martingale problem. Let us explain this in more detail, using the notation from Sections~\ref{sec: approx} and~\ref{sec: subsequence}. For a control rule \(P\) and a function \(w\), set
		\[
		\mathcal{J}_m (P, w) := E^{P} \Big[ w (X_{\tau_{D_m}}) + \int_0^{\tau_{D_m}} \int_\Lambda g (\lambda, X_s) \, M (d \lambda, ds) \Big].
		\]
		For a supersolution \(u\), Lemma~\ref{lem: stochastic lower bound} provides, for every \(m\) with \(x \in D_m\), a {\em single} control rule \(P_m \in \cR_m (x)\) with \(u_m (x) \geq \mathcal{J}_m (P_m, u_m)\), and the passage to the limit only requires that the accumulation points of \((P_m)_{m \geq \J}\) belong to \(\cR (x)\), see Lemma~\ref{lem: tight + closed}. This is a closedness (or upper hemicontinuity) property of the control rules under perturbations of the coefficients.
		In contrast, if \(w_m \in W^2_{p, \on{loc}} (D_{m + 1}) \cap C (D_{m + 1})\) satisfies \(H_m (x, \nabla w_m (x), \nabla^2 w_m (x)) \geq 0\) for a.e. \(x \in D_m\), then Krylov's It\^o formula \cite[Theorem~1.3]{krylov_AOP_21} yields \(w_m (x) \leq \mathcal{J}_m (P, w_m)\) for {\em every} \(P \in \cR_m (x)\). To transfer this inequality to a given \(Q \in \cR (x)\), one has to approximate \(Q\) by control rules of the modified problems, i.e., a lower hemicontinuity property is needed. Such a property may fail if the limiting controlled martingale problem lacks uniqueness.
		For instance, in the notation of Appendix~\ref{ex: non unique}, let \(P^k_n\) be the law of the diffusion with generator \(\frac{1}{2} \on{tr} [ A^k_n \nabla^2 ]\) started at \(x_0\). Then, every accumulation point \(P\) of \((P^k_n)_{n = 1}^\infty\) solves the martingale problem associated with \(\frac{1}{2} \on{tr} [ A \nabla^2 ]\) and satisfies \(E^P [ f (X_{\tau_D}) ] = v^k (x_0)\), cf. Lemmas~\ref{lem: tight + closed} and~\ref{lem: continuity}. As \(v^1 (x_0) \neq v^2 (x_0)\), the sets of accumulation points of \((P^1_n)_{n = 1}^\infty\) and \((P^2_n)_{n = 1}^\infty\) are disjoint. Hence, neither of the two sequences approximates all solutions to the martingale problem associated with \(\frac{1}{2} \on{tr} [ A \nabla^2 ]\).
	\end{remark}

\subsection{Stability with respect to the coefficients} \label{sec: stability}
Based on a stability result for \(L_p\)-viscosity solutions from \cite{KS_22}, for diffusion coefficients satisfying an oscillation condition ensuring uniqueness among \(L_p\)-viscosity solutions, \cite[Proposition~2.9]{C_25_SIMA} established continuity of the value function with respect to perturbations of the coefficients. 
	The following result shows that, without restricting the diffusion coefficients, the value function is always lower semicontinuous with respect to perturbations of the coefficients. As we explain in Remark~\ref{rem: lsc}~(b) below, upper semicontinuity does not hold in general without additional assumptions. 
	
\begin{proposition} \label{coro: lsc}
	Assume that Standing Assumption~\ref{SA: main1} holds.
For \(n \in \mathbb{N}\), let 
\[
a^{(n)} \colon \Lambda \times \bR^d \to \bS, \quad b^{(n)} \colon \Lambda \times \bR^d \to \bR^d
\] 
be Borel functions that satisfy Standing Assumption~\ref{SA: main1} with the same constant \(\delta\) as \((a, b)\) and functions \(\mathfrak{c}_K^{(n)}, w_K^{(n)}\) and \(\b^{(n)}\) that may depend on \(n\).
Moreover, assume that Standing Assumption~\ref{SA: main2} holds with \(\g \in L_{d_0} (D)\) for \(d_0 = d_0 (d, \delta, \| \b \|_{L_d (\bR^d)} \vee \sup_{n \in \mathbb{N}} \| \b^{(n)} \|_{L_d (\bR^d)}) \in (d/2, d)\) as in Theorem~\ref{theo: main1}, which applies for all coefficients simultaneously, that \(\b^{(n)} \to \b\) in \(L_{d} (\bR^d)\), with \(\b \in L_d (\bR^d)\) as in Standing Assumption~\ref{SA: main1}, and
\begin{align} \label{eq: ae convergence coefficients}
	\sup_{\lambda \in \Lambda} \Big( \| a^{(n)} (\lambda, x) - a (\lambda, x) \| + \| b^{(n)} (\lambda, x) - b (\lambda, x) \| \Big) \to 0, \quad n \to \infty, \quad \text{for a.e. } x \in \bR^d.
\end{align}
Denote by \(\cR^{(n)} (x)\) the set of relaxed control rules with initial value \(x\) associated with \((a^{(n)}, b^{(n)})\), and by \(v^{(n)}\) the value function defined as in \eqref{eq: value function D} with \(\cR (x)\) replaced by \(\cR^{(n)} (x)\). Then, it holds that
\[
\liminf_{n \to \infty} v^{(n)} (x) \geq v (x), \quad x \in \oD.
\]
In particular, every pointwise accumulation point of \((v^{(n)})_{n = 1}^\infty\) dominates \(v\).
\end{proposition}

\begin{proof}
We only sketch the argument.
First, \(v^{(n)} (x) = f (x) = v (x)\) for all \(x \in \partial D\) is clear by construction. Take \(x \in D\) and choose a subsequence \((n_j)_{j = 1}^\infty\) such that 
\[
\liminf_{n \to \infty} v^{(n)} (x) = \lim_{j \to \infty} v^{(n_j)} (x).
\] 
By Proposition~\ref{prop: optimal}, applied with \((a^{(n)}, b^{(n)})\) instead of \((a, b)\), there exists a relaxed control rule \(P_n \in \cR^{(n)} (x)\) such that
\[
v^{(n)} (x) = E^{P_n} \Big[ f (X_{\tau_D}) + \int_0^{\tau_D} \int_\Lambda g (\lambda, X_s) \, M (d \lambda, ds) \Big].
\]
Using that the constant \(N = N (d, \delta, \| \b \|_{L_d (\bR^d)} \vee \sup_{n \in \mathbb{N}} \| \b^{(n)} \|_{L_d (\bR^d)}) > 0\) from \cite[Theorem~2.7]{krylov_AOP_21} and \cite[Corollary~1.2]{krylov_21} can be taken uniform for all coefficients simultaneously, we may argue as in the proof of Lemma~\ref{lem: tight + closed} below 
to conclude that \((P_{n_j})_{j = 1}^\infty\) is tight and that all of its accumulation points are in \(\cR (x)\).
The only change in the proof is that the term \(\| (1 + \b ) \1_{G_m} \|_{L_d (\bR^d)}\) in \eqref{eq: long estimate lem 4.4} has to be replaced by \(\| \sup_{\lambda \in \Lambda} | L (\lambda, \, \cdot \,, \nabla \varphi, \nabla^2 \varphi) - L^{(n_j)} (\lambda, \, \cdot \,, \nabla \varphi, \nabla^2 \varphi) | \|_{L_d (\bR^d)}\), with \(L^{(n)}\) defined as in \eqref{eq: L} for the coefficients \((a^{(n)}, b^{(n)})\). The assumptions \(\b^{(n)} \to \b\) in \(L_d (\bR^d)\) and \eqref{eq: ae convergence coefficients}, together with the dominated convergence theorem, ensure that it vanishes. Passing to a further subsequence, we may assume that \(P_{n_j} \to P\) for some \(P \in \cR (x)\). Arguing as in the proof of Lemma~\ref{lem: continuity} below, we obtain that
\[
\liminf_{n \to \infty} v^{(n)} (x) = \lim_{j \to \infty} v^{(n_j)} (x) = E^P \Big[ f (X_{\tau_D}) + \int_0^{\tau_D} \int_\Lambda g (\lambda, X_s) \, M (d \lambda, ds) \Big]  \geq v (x). 
\] 
The final claim is an immediate consequence.
\end{proof} 

\begin{remark} \label{rem: lsc}
(a) Proposition~\ref{coro: lsc} has a direct interpretation for approximation schemes based on regularization. If, for instance, the diffusion coefficient \(a\) is mollified in the space variable, then, in the limit, the value functions of the regularized control problems can only overestimate \(v\). Under the hypotheses of the stability result \cite[Proposition~2.9]{C_25_SIMA}, the sequence \((v^{(n)})_{n = 1}^\infty\) converges, along subsequences, uniformly to an \(L_p\)-viscosity solution of \eqref{eq: PDE}. In this case, Proposition~\ref{coro: lsc} is also a consequence of Theorem~\ref{theo: main1}.

(b) Lower semicontinuity under mollification of the drift and running cost for finite time-horizon control problems with additive noise and Borel action spaces has been established in \cite[Section~3.1]{arXiv:2509.01054}. This setting corresponds to a parabolic HJB equation. By means of a counterexample with a drift that is discontinuous in the control variable, it was shown that upper semicontinuity may fail in general.
For countable action spaces, also upper semicontinuity was proved.
Under the remaining assumptions of that work, the resulting local uniform convergence also holds for compact Polish action spaces when both the drift and running cost are continuous in the control variable, as then one may restrict to a countable dense subset of the action space. In the setting of \cite{arXiv:2509.01054}, the noise is additive and the associated HJB equation admits a unique Sobolev solution by \cite[Theorem~2.1]{arXiv:2509.01054}.

In the present setting, where the diffusion coefficient is merely measurable, uniqueness may fail, and the inequality in Proposition~\ref{coro: lsc} can be strict even under the above assumptions on the action space and continuity in the control variable. Indeed, in the example from Appendix~\ref{ex: non unique}, the value functions associated with the smooth coefficients \((1-\alpha\lambda)A^1_n\) and \((1-\alpha\lambda)A^2_n\) converge, along a subsequence, to two \(L_p\)-viscosity solutions \(w^1\) and \(w^2\) of \eqref{eq: PDE} with \(w^1(x_0)\neq w^2(x_0)\). By Theorem~\ref{theo: main1}, both dominate the value function \(v\), so that at least one of them is strictly larger than \(v\) at \(x_0\). 

It would be interesting to know whether, under the hypotheses of Proposition~\ref{coro: lsc}, \(v\) can always be obtained as a limit of value functions of control problems with drift and
diffusion coefficients that are continuous in the space variable. In this case, \(v\) would be characterized as the smallest such limit.
\end{remark}

\subsection{Concave Hamiltonians} \label{sec: concave}

Before we turn to the proof of our main result in the following section, we explain that many concave uniformly elliptic Hamiltonians have an HJB representation of the form \eqref{eq: Hamiltonian} with coefficients satisfying Standing Assumptions~\ref{SA: main1} and \ref{SA: main2}. The proof of the following proposition can be found in Appendix~\ref{app: pf prop}.

\begin{proposition}  \label{prop: representation} 
	Let \(G \colon \bR^d \times \bR^d \times \mathbb{S}^d \to \bR\) be a map with the following properties:
	\begin{enumerate}
		\item[(i)] For every \(q \in \bR^d\) and \(A \in \mathbb{S}^d\), the map \(x \mapsto G (x, q, A)\) is Borel measurable. 
		\item[(ii)] For every \(x \in \bR^d\), the map \((q, A) \mapsto G (x, q, A)\) is concave and positively homogeneous, i.e., \(G (x, t q, t A) = t \, G (x, q, A)\) for all \(t > 0\).
		\item[(iii)]
		For every \(x \in \mathbb{R}^d\) and \(A \in \mathbb{S}^d_+\), \(G (x, 0, A) \geq 0\).
	\end{enumerate} 
	Furthermore, let \(F \colon \bR^d \times \bR^d \times \mathbb{S}^d \to \bR\) be a map such that, for every \(q \in \bR^d\) and \(A \in \mathbb{S}^d\), the map \(x \mapsto F (x, q, A)\) is Borel measurable, for every \(x \in \bR^d\), the map \((q, A) \mapsto F (x, q, A)\) is concave, and 
	\begin{align} \label{eq: F-G}
		| F (x, q, A) - G (x,q, A) | \leq \mathfrak{g} (x), \quad (x, q, A) \in \bR^d \times \bR^d \times \mathbb{S}^d, 
	\end{align}
	for a Borel function \(\mathfrak{g} \colon \bR^d \to \bR_+\).
	Then, there exists a compact Polish space \(\Gamma\) and measurable functions \((\overline{a}, \overline{b}, \overline{g}) \colon \Gamma \times \bR^d \to \mathbb{S}^d_+ \times \bR^d \times \bR\) such that 
	\[
	F (x, q, A) = \inf_{\gamma \in \Gamma} \Big\{ \langle \overline{b} (\gamma, x), q \rangle + \tfrac{1}{2} \on{tr} \big[ \overline{a} (\gamma, x) A \big] + \overline{g} (\gamma, x) \Big\}, 
	\] 
	and, for a constant \(c > 0\) depending only on \(d\), 
	\begin{equation} \label{eq: for (a) 0}
	\begin{split}
		\| \overline{a} (\gamma, x) - \overline{a}(\gamma', x) \| + \| \overline{b} ( \gamma, x) - \overline{b} (\gamma', x) \| \leq c\,  \sup_{\|A\| + \|q\| \leq 1} | G (x, q, A) | \, d_\Gamma (\gamma, \gamma'), 
	\end{split}
	\end{equation} 
	and 
	\begin{align} \label{eq: properties g}
	|\overline{g} (\gamma, x) | \leq \mathfrak{g} (x), \qquad  | \overline{g} (\gamma, x) - \overline{g} (\gamma', x) |  \leq  c \, \g (x) d_\Gamma (\gamma, \gamma'),
	\end{align} 
	for all \(x \in \bR^d\) and \(\gamma, \gamma' \in \Gamma\).
	Moreover, for any \(x \in \bR^d\), if 
	\begin{align} \label{eq: bound 1}
		| G (x, q, 0) | \leq \b (x)  \| q \|, \quad q \in \mathbb{R}^d, 
	\end{align} 
	then \(\| \overline{b} (\gamma, x) \| \leq \b (x)\), and if
	\begin{align} \label{eq: bound 2}
		\frac{\delta}{2} \on{tr} \big[ S\big] \leq G (x, 0, A + S) - G (x, 0, A) \leq \frac{1}{2\delta} \on{tr} \big[ S \big], \quad A \in \mathbb{S}^d, \, S \in \mathbb{S}^d_+, 
	\end{align} 
	then \(\overline{a} (\gamma, x) \in \mathbb{S}_\delta^d\). Under these two conditions, we also have 
	\begin{align} \label{eq: for (a)}
		\sup_{\|A\| + \|q\| \leq 1} | G (x, q, A) | \leq C (\delta, d) + \b (x), \quad x \in \bR^d,
	\end{align} 
	for a constant \(C\) depending only on \(\delta\) and \(d\).
	In case \(F = G\), we can take \(\mathfrak{g} \equiv 0\), and then \(\overline{g} \equiv 0\).
\end{proposition}

If \eqref{eq: bound 1} holds for a non-negative Borel \(\b \in L_d (\bR^d)\) and \eqref{eq: bound 2} holds with a uniform constant \(\delta \in (0, 1)\) for all \(x \in \bR^d\), then the pair \((\overline{a}, \overline{b})\) has all properties described in Standing Assumption~\ref{SA: main1}, where part~(a) follows from \eqref{eq: for (a) 0} and \eqref{eq: for (a)}. Further, \(\overline{g}\) also satisfies Standing Assumption~\ref{SA: main2}, and, if in addition \(\g \in L_p (D)\) for some \(p \in (d_0, d)\), Theorem~\ref{theo: main1} can be applied with the triplet \((\overline{a}, \overline{b}, \overline{g})\).
To put all this in context, Proposition~\ref{prop: representation} shows that Theorem~\ref{theo: main1} provides a stochastic representation of the minimal \(L_p\)-viscosity supersolution for a wide range of nonlinear Dirichlet problems of the form
\begin{align*}
	\begin{cases} 
		- F (x, \nabla u, \nabla^2 u) = 0, & x \in D, \\ u (x) = f (x), & x \in \partial D, 
	\end{cases} 
\end{align*}
with concave Hamiltonians \(F\) that differ from a concave and positively homogeneous uniformly elliptic Hamiltonian \(G\) with \(L_d\)-drift bound by at most a Borel function \(\g\in L_p (D)\). 

\section{Proof of Theorem~\ref{theo: main1}} \label{sec: proof}

Throughout this section, let \(d_0 = d_0 (d, \delta, \| \b \|_{L_d (\bR^d)}) \in (d/2, d)\) be large enough such that \eqref{eq: krylov estimate} and all relevant results from \cite{C_25_SIMA, KS_22, krylov_CPDE_20, krylov_AOP_21} hold. Further, take \(p \in (d_0, d)\) and assume that \(\g \in L_p (D)\).

\subsection{The viscosity solution property} \label{sec: VP}
In this first subsection, we show that \(v\) is truly an \(L_p\)-viscosity solution to the Dirichlet problem~\eqref{eq: PDE}. We deduce this claim from \cite[Theorem~2.5]{C_25_SIMA}, stability of \(L_p\)-viscosity solutions \cite{KS_22}, and convergence of the value function that follows from a Krylov estimate~\cite{krylov_AOP_21}. 

First, we extend \(g\) and \(\mathfrak{g}\) to the space \(\bR^d\) by setting \(g (\lambda, x) := 0\) and \(\mathfrak{g} (x) := 0\) for \(\lambda \in \Lambda\) and \(x \in \bR^d \setminus D\). We clearly have \(\mathfrak{g} \in L_p (\bR^d)\). Let \(\eta \colon \bR^d \to \bR_+\) be the standard mollification kernel 
\[
\eta (x) := \begin{cases} \displaystyle c \exp \Big\{ \frac{1}{\|x\|^2 - 1} \Big\}, & \|x\| < 1, \\ 0, & \|x\| \geq 1, \end{cases}
\]
where \(c > 0\) is such that \(\int_{\bR^d} \eta (x) \, dx  = 1\). For \(h > 0\), set 
\[
g_h (\lambda, x) := \int_{\{\|z\| < 1\}} \eta (z) g (\lambda, x - hz) \, dz = h^{-d} \int_{D} \eta \Big( \frac{x - z}{h} \Big) g (\lambda, z) \, dz, 
\] 
for \(\lambda \in \Lambda\) and \(x \in \bR^d\). By virtue of Standing Assumption~\ref{SA: main2}, each \(g_h\) is jointly continuous and bounded by 
\[
\sup_{\substack{\lambda \in \Lambda \\ x \in D}} | g_h (\lambda, x) | \leq c h^{-d} \int_{D} \mathfrak{g} (z) \, dz < \infty.
\] Furthermore, 
\begin{align} \label{eq: approx in Lp}
\Big\| \sup_{\lambda \in \Lambda} | g (\lambda, \, \cdot \,) - g_h (\lambda, \, \cdot \, )| \Big\|_{L_p (D)} \to 0, \quad h \searrow 0.
\end{align}
To prove this, take \(\varepsilon > 0\). By the compactness of \(\Lambda\) and Standing Assumption~\ref{SA: main2}, there are \(K \in \mathbb{N}\) points \(\lambda_1, \dots, \lambda_K \in \Lambda\) such that, for all \(x \in D\), 
\[
\sup_{\lambda \in \Lambda} | g (\lambda, x) - g_h (\lambda, x) | \leq \max_{i = 1}^K | g (\lambda_i, x) - g_h (\lambda_i, x) | + \varepsilon\, \Big(\mathfrak{g} (x) + \int_{\{\|z\| \leq 1\}} \eta (z) \mathfrak{g} (x - hz) \, dz \Big).
\] 
Hence, using the triangle and Minkowski's inequalities, we obtain that 
\begin{align*}
	\Big\| \sup_{\lambda \in \Lambda} | g (\lambda, \, \cdot \,) - g_h (\lambda, \, \cdot \, )| \Big\|_{L_p (D)} &\leq \Big( \sum_{i = 1}^K \| g (\lambda_i, \, \cdot \,) - g_h (\lambda_i, \, \cdot \, )\|_{L_p (D)}^p \Big)^{1/p} + 2 \varepsilon \| \mathfrak{g} \|_{L_p (D)}.
\end{align*}
In view of \cite[Lemma~7.2, p. 148]{GT_01}, this estimate implies \eqref{eq: approx in Lp}.
Set 
\[
v_h (x) := \inf_{P \in \cK (x)} E^P \Big[ f (X_{\tau_D}) + \int_0^{\tau_D} \int_{\Lambda} g_h (\lambda, X_s) \, M (d \lambda, ds) \Big], \quad x \in \oD. 
\] 
By \cite[Theorem~2.5]{C_25_SIMA}, the function \(v_h\) is an \(L_p\)-viscosity solution to the nonlinear Dirichlet problem 
\begin{align*}
	\begin{cases}	- H_h (x, \nabla u (x), \nabla^2 u (x)) = 0, & x \in D,  \\
		u (x) = f (x), & x \in \partial D,  \end{cases}
\end{align*}
where 
\[
H_h (x, q, A) :=  \inf_{\lambda \in \Lambda}  \Big\{ \langle b (\lambda, x), q\rangle + \tfrac{1}{2} \on{tr} \big[ a (\lambda, x) A \big] + g_h (\lambda, x)\Big\}.
\] 
Next, we show that \(v_h \to v\) uniformly on \(\oD\) as \(h \searrow 0\).
Using the Krylov estimate \eqref{eq: krylov estimate}, for every \(x \in \overline{D}\), we find that 
\begin{align*}
	| v (x) - v_h (x) | &\leq \sup_{P \in \cK (x)} E^P \Big[ \int_0^{\tau_D} \int_\Lambda \, | g (\lambda, X_s) - g_h (\lambda, X_s) | \, M(d \lambda, ds) \Big] 
	\\&\leq N \Big\| \sup_{\lambda \in \Lambda} | g (\lambda, \, \cdot \,) - g_h (\lambda, \, \cdot \, )| \Big\|_{L_p (D)}, 
\end{align*} 
where \(N = N(d, \delta, \|\b\|_{L_d (\bR^d)}, \on{diam} (D)) > 0\). Together with \eqref{eq: approx in Lp}, this estimate implies that \(v_h \to v\) uniformly on \(\oD\) as \(h \searrow 0\).  In particular, \(v \in C (\oD)\), since each \(v_h\) is continuous on \(\oD\) as an \(L_p\)-viscosity solution (see Definition~\ref{def: VS elliptic}).
Finally, using \eqref{eq: approx in Lp} again, we conclude from the stability result \cite[Proposition~4.6]{KS_22} that \(v\) is an \(L_p\)-viscosity solution to \eqref{eq: PDE}. 

In the remainder of this section, we prove that it is the minimal \(L_p\)-viscosity supersolution. This proof is split into several steps, given in the following subsections.

\subsection{Approximate solutions} \label{sec: approx}
Let \((D_n)_{n = 1}^\infty\) be a sequence of bounded domains with smooth boundaries such that \(D_n \subset \subset D_{n + 1} \subset \subset D\) and \(D = \bigcup_{n = 1}^\infty D_n\), see \cite[Proposition~8.2.1]{D_08}. 
\begin{SA} \label{SA: p, u}
For the rest of this proof, let \(u\) be an arbitrary \(L_p\)-viscosity supersolution to the Dirichlet problem \eqref{eq: PDE}. 
\end{SA}
\begin{lemma} \label{lem: main approx}
	Let \(u\) be as in Standing Assumption~\ref{SA: p, u}.
		There exists a sequence \(H_m \colon \bR^d \times \bR^d \times \mathbb{S}^d \to \bR\) of HJB Hamiltonians with the following properties:
	\begin{enumerate}
		\item For all \(x, q \in \bR^d\) and \(A \in \mathbb{S}^d\), 
		\[
		H_m (x, q, A) = \inf_{\lambda \in \Lambda} \Big\{ \langle b_m (\lambda, x), q \rangle + \tfrac{1}{2} \on{tr} \big[ a_m (\lambda, x) A \big] + g (\lambda, x)\Big\}, 
		\] 
		where \(a_m \colon \Lambda \times \bR^d \to \mathbb{S}_\delta^d\) and \(b_m \colon \Lambda \times \bR^d \to \bR^d\) are Borel functions satisfying the following properties: 
		\begin{enumerate}
		\item for every compact set \(K \subset \bR^d\), 
		 		\[
		 \| a_m (\lambda, x) - a_m (\lambda', x) \| + \| b_m (\lambda, x) - b_m (\lambda', x) \| \leq \mathfrak{c}_K (x) \, w_K \big(d_\Lambda (\lambda, \lambda') \big) 
		 \] 
		 for all \(\lambda, \lambda' \in \Lambda\) and \(x \in K\), where \(\mathfrak{c}_K\) is the function from Standing Assumption~\ref{SA: main1}~(a);
		\item \(\| b_m (\lambda, x) \| \leq \b (x)\) for all \(m \in \mathbb{N}, \lambda \in \Lambda\) and \(x \in \bR^d\); 
		\end{enumerate} 
		\item There are measurable sets \((G_m)_{m = 1}^\infty\) with \(G_m \subset D_m\), \(\lim_{m \to \infty} \on{meas} \, (G_m) = 0\), and 
		\[
		b_m (\lambda, x) = b (\lambda, x), \qquad a_m (\lambda, x) = a (\lambda, x),
		\] 
		for all \(\lambda \in \Lambda\) and \(x \not \in G_m\); 
		\item There are functions \(u_m \in  W^2_{p, \on{loc}} (D_{m + 1}) \cap C (D_{m + 1})\) such that
		\[
		H_m (x, \nabla u_m (x), \nabla^2 u_m (x)) \leq 0 \quad \text{for a.e. \(x \in D_m\)}; 
		\] 
		\item \(u_m \leq u\) on \(D_{m + 1}\) and, as \(m \to \infty\), 
		\[
		\sup_{y \in \overline{D}_m} | u_m (y) - u (y) | \to 0. 
		\] 
	\end{enumerate} 
\end{lemma}

\begin{proof}
	The proof is split into three steps. In the first we construct the approximation sequence \((u_m)_{m = 1}^\infty\) and in the second we construct the associated Hamiltonians \((H_m)_{m = 1}^\infty\). In the third and last step, we put the pieces together. 
	
	\smallskip 
	{\em Step 1:}
	We start by introducing some standard notation. Namely, we define the Pucci extremal operators by 
	\begin{align*}
		\mathcal{P}^+ (A) := - \frac{\delta}{2} \on{tr} \big[ A^+ \big] + \frac{1}{2\delta} \on{tr} \big[ A^- \big], \qquad
		\mathcal{P}^- (A) := - \frac{1}{2 \delta} \on{tr} \big[ A^+ \big] + \frac{\delta}{2} \on{tr} \big[ A^- \big],  
	\end{align*}
	where \(A^+\) and \(A^-\) are the positive and negative parts of \(A\), respectively. Notice that, as in \cite{CCKS_96,KS_22}, \(\mathcal{P}^\pm\) are decreasing in \(A\).
	Furthermore, we set 
	\[
	\beta (x) := \sup_{\lambda \in \Lambda} \| b (\lambda, x) \|, \qquad \ug (x) := \inf_{\lambda \in \Lambda} g (\lambda, x), 
	\] 
	which are Borel functions by virtue of Standing Assumptions~\ref{SA: main1}~(a) and \ref{SA: main2}, and the compactness of the action space \(\Lambda\). 
	
		Since \(\mathcal{P}^+ (A) = \sup_{a' \in \mathbb{S}^d_\delta} \, ( - \frac{1}{2} \on{tr} [ a' A ] )\) and \(a (\lambda, x) \in \mathbb{S}^d_\delta\), we have
	\begin{align} \label{eq: H Pucci}
		- H (x, q, A) = \sup_{\lambda \in \Lambda} \Big\{ - \langle b (\lambda, x), q \rangle - \tfrac{1}{2} \on{tr} \big[ a (\lambda, x) A \big] - g (\lambda, x) \Big\} \leq \mathcal{P}^+ (A) + \beta (x) \| q \| - \ug (x)
	\end{align}
	for all \(x \in D\), \(q \in \bR^d\) and \(A \in \mathbb{S}^d\). In particular, \(u\) is also an \(L_p\)-viscosity supersolution of \(\mathcal{P}^+ (\nabla^2 u) + \beta \| \nabla u \| - \ug = 0\) in \(D\), i.e., whenever \(\varphi \in W^2_{p, \textup{loc}} (D)\) is such that \(u - \varphi\) has a local minimum at \(x_0 \in D\), then
	\begin{align} \label{eq: Pucci super}
		\lim\nolimits_{r \searrow 0} \on{ess \, sup}_{x \in B_r (x_0)} \big( \mathcal{P}^+ (\nabla^2 \varphi (x)) + \beta (x) \| \nabla \varphi (x) \| - \ug (x) \big) \geq 0.
	\end{align}
	
	For \(k, n, m \in \mathbb{N}\), in view of \cite[Theorem~1.2]{krylov_CPDE_20}, there exists a function \(v_{k, n}^m \in W^2_{p, \on{loc}} (D_m) \cap C (\oD_m)\) such that 
	\begin{align} \label{eq: def HF}
	\begin{cases} \mathcal{P}^+ (\nabla^2 v_{k, n}^m (x) ) + \beta (x) \| \nabla v_{k, n}^m (x) \| + n (v_{k, n}^m (x) - u (x))_+ = \ug (x) + k& \text{a.e. in \(D_m\)},
		\\ v^m_{k, n} = u &\text{on } \partial D_m, 
		\end{cases} 
	\end{align} 
	where \(x_+ := \max \{x, 0\}\).
	To clarify, we used here \cite[Theorem~1.2]{krylov_CPDE_20} to produce the solution \(- v^m_{k, n}\), whose corresponding equation satisfies the assumptions from \cite{krylov_CPDE_20}.
	
	Before going into the details, let us explain the idea behind the construction. In Step~1, we let \(n \to \infty\) for fixed \(k\) and \(m\). Where \(v^m_{k,n} > u + k/n\), the penalization \(n (v^m_{k,n} - u)_+\) makes \(v^m_{k,n}\) a strict subsolution of the Pucci equation, of which \(u\) is a supersolution by \eqref{eq: Pucci super}. Applying the definition of an \(L_p\)-viscosity supersolution to \(u\), with \(v^m_{k,n}\) as a test function at a maximum point of \(v^m_{k,n} - u\), shows that this region is empty, i.e., \(v^m_{k,n} \leq u + k/n\). We then show that, along a subsequence, \(v^m_{k,n}\) converges locally uniformly to a function \(v^m_k \in W^2_{p,\textup{loc}}(D_m) \cap C(D_m)\). In the limit, \(v^m_k \leq u\) and, as the penalization vanishes where \(v^m_k < u\), we obtain \(\mathcal{P}^+(\nabla^2 v^m_k) + \beta \|\nabla v^m_k\| = \ug + k\) a.e. on the non-contact set \(\{v^m_k < u\}\). The large source term \(k\) pushes \(v^m_k\) towards \(u\), so that, up to an error of order \(1/m\) coming from a smooth approximation of \(u\), \(v^m_k\) approaches \(u\) as \(k \to \infty\). In Step~2, we observe that \(v^m_k\) inherits the supersolution property of \(u\) for the original Hamiltonian on the contact set \(\{v^m_k = u\}\), so that the Hamiltonian only has to be modified on the non-contact set. 
	To see that the measure of this set is small, note that, a.e. in \(D_m\),
	\[
	\1_{\{ v^m_k < u \}} = \frac{1}{k} \big( \mathcal{P}^+(\nabla^2 v^m_k) + \beta \|\nabla v^m_k\| - \ug \big) \1_{\{ v^m_k < u \}} \leq \frac{1}{k} \big( \mathcal{P}^+(\nabla^2 v^m_k) + \beta \|\nabla v^m_k\| - \ug \big).
	\]
	We integrate this inequality along the diffusion whose coefficients realize the Pucci operator along \(v^m_k\). As \(v^m_k\) is bounded uniformly in \(k\), It\^o's formula bounds the time integral of the right-hand side by \(C/k\), with \(C\) independent of \(k\). Hence, the diffusion spends an expected time of order \(1/k\) in the non-contact set, and a Krylov--Safonov type estimate shows that its Lebesgue measure vanishes as \(k \to \infty\).
		
We now make this program rigorous. We claim that \(v^m_{k, n} \leq u + k / n\) on \(\oD_m\). To prove this, assume for contradiction that
\[
\max_{\oD_m} \, (v^m_{k, n} - u) > \frac{k}{n}.
\]
As \(v^m_{k, n} = u\) on \(\partial D_m\), there exists an \(x_0 \in D_m\) such that \(v^m_{k, n} (x_0) - u (x_0) = \max_{\oD_m} \, (v^m_{k, n} - u)\) and, by continuity, there exist a radius \(r_0 > 0\) with \(\overline{B}_{r_0} (x_0) \subset D_m\) and a number \(\eta > 0\) such that
\[
v^m_{k, n} - u \geq \frac{k}{n} + \eta \text{ on } B_{r_0} (x_0).
\]
As a consequence, using \eqref{eq: def HF}, a.e. on \(B_{r_0} (x_0)\),
\begin{align} \label{eq: from 4.4}
	\mathcal{P}^+ (\nabla^2 v_{k, n}^m ) + \beta \| \nabla v_{k, n}^m \| - \ug = k - n ( v^m_{k, n} - u )_+ \leq k - (k + n \eta) = - n \eta.
\end{align}
Take a cutoff function \(\chi \in C^\infty_c (D_m)\) with \(\chi = 1\) on \(B_{r_0} (x_0)\), and set \(\varphi := \chi v^m_{k, n}\) on \(D_m\) and \(\varphi := 0\) on \(D \setminus D_m\). Then, \(\varphi \in W^2_{p, \textup{loc}} (D)\), \(\varphi = v^m_{k, n}\) on \(B_{r_0} (x_0)\), and \(u - \varphi\) has a local minimum at \(x_0\). By \eqref{eq: from 4.4}, the left-hand side of \eqref{eq: Pucci super} is at most \(- n \eta < 0\), which is a contradiction. We conclude that \(v^m_{k, n} \leq u + k / n\) on \(\oD_m\).
	In particular, 
	\begin{align} \label{eq: new inequalities}
	\ug \leq 	\mathcal{P}^+ (\nabla^2 v_{k, n}^m ) + \beta (x) \| \nabla v_{k, n}^m \| \leq \ug + k \quad \text{a.e. on \(D_m\)}.
	\end{align} 
	
	Next, we derive a lower bound for \(v^m_{k, n}\).
	Using \cite[Theorem~2.1]{krylov_CPDE_20} with \(- v^m_{k, n}\), we obtain that 
	\begin{equation} \label{eq: lower bound} \begin{split}
	- v^m_{k, n} & \leq N \| (- \mathcal{P}^- (\nabla^2 (- v^m_{k, n}  )) + \beta \| \nabla (- v^m_{k, n} ) \|)_- \|_{L_p (D_m)} + \sup_{\partial D_m} (- v^m_{k, n})_+ 
	\\&= N \| (\mathcal{P}^+ (\nabla^2 v^m_{k, n} ) + \beta \| \nabla v^m_{k, n} \|)_- \|_{L_p (D_m)} + \sup_{\partial D_m} (- v^m_{k, n})_+ 
	\\&\leq N \| \ug_- \|_{L_p (D)} + \| u \|_\infty,
	\end{split}
	\end{equation} 
	where \(x_- := \max \{- x, 0\}\).
	We conclude that \(v^m_{k, n}\) is bounded on \(D_m\) by a constant independent of~\(n\). 
	Using this observation, we get from \cite[Theorem~1.1]{krylov_CPDE_20}, applied to \(- v^m_{k, n}\), that, for every open domain \(D_m' \subset \subset D_m\), the sequence \( (v^m_{k, n})_{n = 1}^\infty\) is bounded in \(W^2_{p} (D_m')\). To see this, notice first that, on \(\oD_m\),  
	\[
	\big| n ( v^m_{k, n}- u )_+ - \ug - k \big| \leq n ( v^m_{k, n} - u )_+ + | \ug | + k \leq n \Big( u + \frac{k}{n} - u \Big )_+ + |\ug| + k \leq |\ug| + 2k.
	\] 
	Hence, using the notation of \cite{krylov_CPDE_20}, we may take
	\[
	G (\mathsf{u}, x) = \beta (x) \| \nabla \mathsf{u} \| + n ( v^m_{k, n} (x) - u (x) )_+ - \ug (x) - k, \quad K_0 = 0, \quad \bar{G} (x) = | \ug (x) | + 2 k,
	\] 
	to obtain from \cite[Theorem~1.1]{krylov_CPDE_20} that
	\begin{align*}
		\| v^m_{k, n} \|_{W^2_p (D_m')} \leq C \Big( \| |\ug| + 2k \|_{L_p (D_m)} + \sup_{D_m} |v^m_{k, n} | \Big), 
	\end{align*}
	where \(C > 0\) is independent of \(n\). 
	Hence, as \(p > d / 2\), the Sobolev embedding theorem yields that there exists a function \(v^m_k \in W^2_{p, \on{loc}} (D_m) \cap C (D_m)\) such that, along a subsequence, 
	\(
	v^m_{k, n_j} \to v^m_k
	\) 
	locally uniformly in \(C (D_m)\). In particular, 
	\begin{align} \label{eq: order} 
		v^m_k \leq u.
	\end{align}  
	
	Recalling \eqref{eq: new inequalities}, we deduce from \cite[Propositions~4.5, 4.6, 5.3]{KS_22} that also 
	\begin{align} \label{eq: a.e. super}
\ug \leq 	\mathcal{P}^+ (\nabla^2 v_{k}^m ) + \beta (x) \| \nabla v_{k}^m \| \leq \ug + k \quad \text{a.e. on \(D_m\)}.
	\end{align} 
	Define the open set \[G_k^m := \{ v^m_k < u \}.\]
	For every relatively compact set \(U \subset \subset G_k^m\), we have \(\inf_U (u - v^m_k) > 0\). Hence, for sufficiently large~\(j\), as \(v^m_{k, n_j} \to v^m_k\) uniformly on \(\overline{U}\), we get \((v^m_{k, n_j} - u)_+ = 0\) on \(U\).
	Using \eqref{eq: def HF} and \cite[Propositions~4.5, 4.6, 5.3]{KS_22}, we conclude that
	\[
		\mathcal{P}^+ (\nabla^2 v^m_k) + \beta \| \nabla v^m_k \| = \ug + k \quad \text{a.e. on \(U\)}. 
	\]
	This entails that  
	\begin{align} \label{eq: eq limit}
	\mathcal{P}^+ (\nabla^2 v^m_k) + \beta \| \nabla v^m_k \| = \ug + k \quad \text{a.e. on \(G^m_k\)}. 
	\end{align} 
	
	Take a function \(\varphi^m \in C^\infty (D_{m + 1})\) such that 
	\[
	u - \frac{1}{m} \leq \varphi^m \leq u \text{ on } \oD_m.
	\] 
	The existence of such a function is guaranteed by a standard mollification result \cite[Lemma~7.1, p. 147]{GT_01}. Namely, there exists a function \(\psi^m \in C^\infty (D_{m + 1})\) such that \(| \psi^m - u | \leq \frac{1}{2m}\) on \(\oD_m\). Then, take \(\varphi^m := \psi^m - \frac{1}{2m}\), which satisfies
	\[
	u - \frac{1}{m} \leq u - \frac{1}{2m} - | \psi^m - u | \leq \varphi^m \leq \psi^m - | \psi^m - u | \leq u. 
	\] 
	
	Recalling that \(v^m_k \leq u\) and \(u - \varphi^m \leq \frac{1}{m}\), we get that 
	\[
	| u - v^m_k | = u - v^m_k \leq \frac{1}{m} + \varphi^m - v^m_k.
	\] 
	To estimate the last term, set \(U := \{ x \in D_m \colon v^m_{k, n} (x) < \varphi^m (x) \}\) and \(w := \varphi^m - v^m_{k, n}\). As \(\varphi^m \leq u\), we have \((v^m_{k, n} - u)_+ = 0\) on \(U\), so that \(\mathcal{P}^+ (\nabla^2 v^m_{k, n}) + \beta \| \nabla v^m_{k, n} \| = \ug + k\) a.e. on \(U\) by \eqref{eq: def HF}. Using that \(\mathcal{P}^+ (A + B) \geq \mathcal{P}^+ (A) + \mathcal{P}^- (B)\) for \(A, B \in \mathbb{S}^d\) and \(\| \nabla \varphi^m \| \geq \| \nabla v^m_{k, n} \| - \| \nabla w \|\), we obtain that, a.e. on~\(U\),
	\[
	- \mathcal{P}^- (\nabla^2 w) + \beta \| \nabla w \| \geq - \big( \mathcal{P}^+ (\nabla^2 \varphi^m) + \beta \| \nabla \varphi^m \| - \ug - k \big).
	\]
	Moreover, \(w = 0\) on \(\partial U \cap D_m\) and \(w \leq 0\) on \(\partial D_m\), as \(v^m_{k, n} = u \geq \varphi^m\) there. Hence, the maximum principle \cite[Theorem~2.1]{krylov_CPDE_20}, applied with \(\varphi^m - v^m_{k, n}\) on subdomains of \(U\) as in \eqref{eq: lower bound}, and the fact that \(w \leq 0\) on \(D_m \setminus U\) yield
	\begin{align*}
	\sup_{D_m} \,	(\varphi^m - v^m_{k, n})_+
		\leq N \| (\mathcal{P}^+ (\nabla^2 \varphi^m) +  \beta \| \nabla \varphi^m \| - \ug - k)_+ \|_{L_p (D_m)}.
	\end{align*}
	In summary, with \(v^m_{k, n_j} \to v^m_k\), we get that 
	\begin{align} \label{eq: for convergence}
		| u - v^m_k | \leq \frac{1}{m} + N \| (\mathcal{P}^+ (\nabla^2 \varphi^m) + \beta \| \nabla \varphi^m \| - \ug - k)_+ \|_{L_p (D_m)}.
	\end{align} 
	As \(\beta \in L_d (D)\) and \(\ug \in L_p (D)\), the dominated convergence theorem yields that 
	\begin{align} \label{eq: for convergence 2}
\| (\mathcal{P}^+ (\nabla^2 \varphi^m) + \beta \| \nabla \varphi^m \| - \ug - k)_+ \|_{L_p (D_m)} \to 0, \quad k \to \infty. 
	\end{align} 

	\smallskip
	{\em Step 2:} 
	Set 
	\[
	\overline{H}^m_k (x, q, A) := \begin{cases} H (x, q, A), & x \not \in G^m_k, \\ - \mathcal{P}^+ (A) - \beta \| q \| + \ug (x), & x \in G^m_k, \end{cases} 
	\] 
	and notice that 
	\begin{align} \label{eq: bound hamiltonians}
	\overline{H}^m_k (x, q, A) \leq H (x, q, A)
	\end{align}
	for all \(x \in D_m, q \in \mathbb{R}^d\) and \(A \in \mathbb{S}^d\). 
	We claim that \(v^m_k\) is an \(L_p\)-viscosity supersolution to 
	\[
	- \overline{H}^m_k (x, \nabla v^m_k, \nabla^2 v^m_k ) = 0.
	\] 
	On the set \(G^m_k\), this follows from the first inequality in \eqref{eq: a.e. super} and \cite[Proposition~4.5]{KS_22}. On the complement of \(G^m_k\), which equals \(\{ v^m_k = u \}\), we can argue as in the proof of \cite[Theorem~4.1]{Jen_Sw_05}.
	Namely, any lower test function for \(v^m_k\) is also a lower test function for \(u\) that coincides with it at the contact point. Thus, the supersolution property of \(u\) propagates to \(v^m_k\) thanks to the inequality \eqref{eq: bound hamiltonians}. Now, \cite[Proposition~5.3]{KS_22} entails that 
	\begin{align} \label{eq: property (2)}
	\overline{H}^m_k (x, \nabla v^m_k (x), \nabla^2 v^m_k (x)) \leq 0 \quad \text{a.e. \(x \in D_m\)}.
	\end{align} 
	
		Denote by \(\Pi^+ (A)\) the spectral projection on the non-negative eigenspace of \(A \in \mathbb{S}^d\), and set, for \(x \in D_m\), 
	\begin{align*}
		a^m_k (x) &:= \delta \Pi^+ (\nabla^2 v^m_k (x)) + \frac{1}{\delta} (\on{id}_d - \Pi^+ (\nabla^2 v^m_k (x))), 
		\\
		b^m_k (x) &:= \begin{cases} - \frac{\beta (x)}{\| \nabla v^m_k (x) \|} \nabla v^m_k (x), & \nabla v^m_k (x) \not = 0,  \\ 0, & \nabla v^m_k (x) = 0, \end{cases} 
	\end{align*}
	where we use Borel representatives of the weak derivatives \(\nabla v^m_k\) and \(\nabla^2 v^m_k\).
	These are Borel maps such that \(\| b^m_k (x) \| \leq \beta (x) \leq \b (x)\) and \(a^m_k (x) \in \mathbb{S}_\delta^d\). Furthermore, for all \(x \in D_m\), 
	\begin{align*}
	-\langle b^m_k (x), \nabla v^m_k (x) \rangle - \tfrac{1}{2} \on{tr} \big[ a^m_k (x) \nabla^2 v^m_k (x) \big] = \mathcal{P}^+ (\nabla^2 v^m_k (x)) + \beta (x) \| \nabla v^m_k (x) \|.
	\end{align*}
	Now, for \(x \in \bR^d, q \in \bR^d\) and \(A \in \mathbb{S}^d\), set 
	\begin{align} \label{eq: definition Hm}
	H^m_k (x, q, A) := \begin{cases} H (x, q, A), & x \not \in G^{m + 1}_k \cap D_m, \\ \langle b^{m + 1}_k (x), q\rangle + \tfrac{1}{2} \on{tr} [ a^{m + 1}_k (x) A ] + \ug(x), & x \in G^{m + 1}_k \cap D_m, \end{cases} 
	\end{align} 
	where it is useful to notice that 
	\begin{equation} \label{eq: inf identity}
		\begin{split}
	\langle b^{m + 1}_k (x), q\rangle &+ \tfrac{1}{2} \on{tr} [ a^{m + 1}_k (x) A ] +\ug(x) 
	\\&= \inf_{\lambda \in \Lambda} \Big\{ \langle b^{m + 1}_k (x), q\rangle + \tfrac{1}{2} \on{tr} [ a^{m + 1}_k (x) A ] + g (\lambda, x) \Big\}, \quad x \in G^{m+ 1}_k \cap D_m.
	\end{split}
	\end{equation}
	From \eqref{eq: property (2)}, we find that 
	\begin{align} \label{eq: final (2)}
		H^m_k (x, \nabla v^{m + 1}_k (x), \nabla^2 v^{m + 1}_k (x)) \leq 0 \quad \text{a.e. \(x \in D_m\)}. 
	\end{align} 
	
	Extend the coefficients \(b^m_k\) and \(a^m_k\) measurably to \(\bR^d\) such that \(\|b^m_k (x)\| \leq \b (x)\) and \(a^m_k (x) \in \mathbb{S}_\delta^d\) for all \(x \in \bR^d\). 	Let \(y_0 \in D_m\), and let \(B = B_r (z_0) \subset\subset D_m\) be a ball such that \(y_0 \in B_{r / 2} (z_0)\). 
	By \cite[Theorem~1.1]{krylov_AOP_21}, there exists a weak solution process \(Y\) to the SDE 
	\[
	d Y_t = b^m_k (Y_t) \, dt + \sqrt{a^m_k (Y_t)} \, d W_t, \quad Y_0 = y_0. 
	\] 
Using \eqref{eq: a.e. super}, \eqref{eq: eq limit} and the Krylov estimate \cite[Theorem~2.8]{krylov_AOP_21}, we obtain 
	\begin{align*}
		E \Big[ \int_0^{\tau_{B} (Y)} \1_{G^m_k} (Y_s)\, ds \Big] \leq \frac{1}{k} E \Big[ \int_0^{\tau_{B} (Y)} \big( \mathcal{P}^+ (\nabla^2 v^m_k (Y_s)) + \beta (Y_s) \| \nabla v^m_k (Y_s) \| - \ug (Y_s) \big) \, ds \Big]. 
	\end{align*}
	Using Krylov's It\^o formula \cite[Theorem~1.3]{krylov_AOP_21}, applied to a function from \(W^2_{p} (\bR^d)\) that coincides with \(v^m_k\) on \(\overline{B}\), the lower bound \eqref{eq: lower bound}, the fact that \(v^m_k \leq u\) from \eqref{eq: order}, and the Krylov estimate \cite[Theorem~2.8]{krylov_AOP_21}, for every \(t \in \bR_+\), 
	we further get that 
	\begin{align*}
		E \Big[ \int_0^{t \wedge \tau_{B} (Y)} \big( \mathcal{P}^+ (\nabla^2 v^m_k (Y_s)) &+ \beta (Y_s) \| \nabla v^m_k (Y_s) \| - \ug (Y_s) \big) \, ds \Big] 
		\\&= E \Big[ v^m_k (y_0) - v^m_k (Y_{t \wedge \tau_{B} (Y)}) - \int_0^{t \wedge \tau_{B} (Y)} \ug (Y_s) \, ds \Big]
		\\&\leq \| u \|_\infty + N \| \ug_- \|_{L_p (D)} + \| u \|_\infty  + E \Big[ \int_0^{\tau_{B} (Y)} \ug_- (Y_s) \, ds \Big]
		\\&\leq 2 \| u \|_\infty + 2N \| \ug_-\|_{L_p (D)}. \phantom \int 
	\end{align*}
	Notice that a.s. the Lebesgue measure of the time the process \(Y\) spends in a Lebesgue null set is zero. This follows from the Krylov estimate \cite[Theorem~2.8]{krylov_AOP_21}. Thus, by virtue of \eqref{eq: a.e. super}, the process 
	\[
	t \mapsto \int_0^{t \wedge \tau_B (Y)}  \big( \mathcal{P}^+ (\nabla^2 v^m_k (Y_s)) + \beta (Y_s) \| \nabla v^m_k (Y_s) \| - \ug (Y_s) \big) \, ds
	\] 
	is a.s. increasing. Hence, the monotone convergence theorem yields that 
	\begin{align*}
	\lim_{t \to \infty} E \Big[ \int_0^{t \wedge \tau_{B} (Y)} \big( &\mathcal{P}^+ (\nabla^2 v^m_k (Y_s)) + \beta (Y_s) \| \nabla v^m_k (Y_s) \| - \ug (Y_s) \big) \, ds \Big] 
	\\&= E \Big[ \int_0^{\tau_{B} (Y)} \big( \mathcal{P}^+ (\nabla^2 v^m_k (Y_s)) + \beta (Y_s) \| \nabla v^m_k (Y_s) \| - \ug (Y_s) \big) \, ds \Big].
	\end{align*} 
	In summary, we showed that 
	\begin{align*}
		E \Big[ \int_0^{\tau_{B} (Y)} \1_{G^m_k} (Y_s)\, ds \Big] \leq \frac{1}{k} \Big( 2 \| u \|_\infty + 2N \| \ug_- \|_{L_p (D)} \Big). 
	\end{align*}
	Thanks to \cite[Theorem~4.1]{krylov_AOP_21}, there are constants \(\mu \geq 1, N > 0\), independent of \(k, B\) and \(y_0\), such that 
	\begin{align} \label{eq: bound measure}
	\Big( \frac{\on{meas} \big (G^m_k \cap B \big)}{\on{meas} \, (B)} \Big)^\mu \leq \frac{N}{r^2} E \Big[ \int_0^{\tau_{B} (Y)} \1_{G^m_k} (Y_s)\, ds \Big] \leq \frac{N}{r^2 k}, 
	\end{align} 
	where the constant \(N\) may change from line to line. 
	
	{\em Step 3:} Using \eqref{eq: for convergence}, \eqref{eq: for convergence 2} and \eqref{eq: bound measure}, covering the compact set \(\oD_m\) by a finite number of balls strictly contained in \(D_{m + 1}\), we may take \(k (m) \in \mathbb{N}\) such that 
	\begin{align} \label{eq: for (4)}
	\on{meas} \big (G^{m + 1}_{k (m)} \cap D_m \big) \leq \frac{1}{m}, \quad \sup_{\oD_m} | v^{m + 1}_{k (m)} - u | \leq \frac{2}{m + 1}. 
	\end{align} 
	Now, set 
	\[
	u_m := v^{m + 1}_{k (m)}, \quad H_m := H^{m}_{k (m)}.
	\] 
	In the previous steps, we proved that (1) -- (4) hold. 
	Namely, for (1) see \eqref{eq: definition Hm} and \eqref{eq: inf identity}, and notice that it holds with 
	\begin{align*}
	b_m (\lambda, x) &:= \begin{cases} b (\lambda, x), & x \not \in G^{m + 1}_{k (m)} \cap D_m, \\ b^{m + 1}_{k (m)} (x), & \text{otherwise}, \end{cases} \\ 
	a_m (\lambda, x) &:= \begin{cases} a (\lambda, x), & x \not \in G^{m + 1}_{k (m)} \cap D_m, \\ a^{m + 1}_{k (m)} (x), & \text{otherwise}. \end{cases}	
	\end{align*} 
	The item (2) holds with \(G_m := G^{m + 1}_{k (m)} \cap D_m\). 
	Finally, for (3) see \eqref{eq: final (2)}; and for (4) see \eqref{eq: order} and \eqref{eq: for (4)}.
	The proof is complete.
\end{proof}

\begin{remark} \label{rem: comparison}
		To relate the construction from Lemma~\ref{lem: main approx} to the literature, Pucci extremal equations were already used in \cite[Section~3]{J_96} and Step~1 of the proof of \cite[Theorem~3.1]{JKS_02}. The term \(n (v_{k, n}^m - u)_+\) is in the spirit of classical approximations for obstacle problems as in \cite[Chapter~4, Section~1.1]{zbMATH03751685} and it plays a role similar to that of the obstacle in the proof of \cite[Theorem~4]{zbMATH05032612} and \cite[Theorem~4.1]{Jen_Sw_05}, and Step~2 in the proof of Lemma~\ref{lem: main approx} adapts the contact-set argument from the proof of the latter. A related use of sets of vanishing measure appears in the proofs of \cite[Theorems~3.30, 3.38]{J_96} and \cite[Theorem~3.1]{JKS_02}, where the approximation involves solution-dependent spatial transformations of the original operators and convex combinations. Here, we do not use such spatial transformations.
\end{remark}

The next lemma relates the approximation sequence to a control representation. 
Let \(\cR_m (x)\) be defined analogously to \(\cR (x)\) with the coefficients \((a, b)\) replaced by \((a_m, b_m)\) from Lemma~\ref{lem: main approx}.
By Lemma~\ref{lem: main approx}~(1), the coefficients \((a_m, b_m)\) satisfy Standing Assumption~\ref{SA: main1} with the same constant \(\delta\) and the same function \(\b\) as \((a, b)\). Consequently, the Krylov estimate \eqref{eq: krylov estimate} holds for all \(P \in \cR_m (x)\), \(x \in \bR^d\) and \(m \in \mathbb{N}\), with the same constants \(d_0\) and \(N\). In particular, the expectation in the following lemma is well-defined.
\begin{lemma} \label{lem: stochastic lower bound}
	Take \(m \in \mathbb{N}\) and \(x \in D_m\). There exists a measure \(P_m \in \cR_m (x)\) such that 
	\begin{align*} 
	E^{P_m} \Big[ u_m (X_{\tau_{D_m}}) + \int_0^{\tau_{D_m}} \int_\Lambda g (\lambda, X_s) \, M (d \lambda, ds) \Big] \leq u_m (x).
	\end{align*} 
\end{lemma}

\begin{proof}
	By the measurable maximum theorem \cite[Theorem~18.19]{charalambos2013infinite}, there exists a Borel function \(\lambda^*_m \colon \bR^d \to \Lambda\) such that, for all \(x \in D_m\),  
	\begin{align*}
	\langle b_m (\lambda^*_m (x), x), \nabla u_m (x) \rangle + \tfrac{1}{2} \on{tr} \big[ a_m (\lambda^*_m (x), x) \nabla^2 u_m (x) \big] &+ g (\lambda^*_m (x), x) 
	\\&= H_m (x, \nabla u_m (x), \nabla^2 u_m (x) ), 
	\end{align*}
	where we use Borel representatives of the weak derivatives \(\nabla u_m\) and \(\nabla^2 u_m\). 
	By \cite[Theorem~1.1]{krylov_AOP_21}, there exists a weak solution process \(Y\) to the SDE 
	\[
	d Y_t = b_m (\lambda^*_m(Y_t), Y_t) \, dt + \sqrt{a_m (\lambda^*_m(Y_t), Y_t)} \, d W_t, \quad Y_0 = x. 
	\] 
	Set \(P := \text{law of } (Y, \delta_{\lambda^*_m(Y_t)} \, dt)\), seen as a probability measure on \((\Omega \times \mathbb{M}, \mathcal{F} \otimes \mathcal{M})\), and notice that \(P \in \cR_m (x)\). 
	Using Krylov's It\^o formula \cite[Theorem~1.3]{krylov_AOP_21}, applied to a function from \(W^2_p (\bR^d)\) that coincides with \(u_m \in W^2_{p, \on{loc}} (D_{m + 1})\) on \(\oD_m \subsetneq D_{m + 1}\), we get, for every \(t \in \bR_+\), that 
	\begin{equation} \label{eq: main inequality main step}
		\begin{split} 
		E^P \Big[ u_m (&X_{t \wedge \tau_{D_m}}) + \int_0^{t \wedge \tau_{D_m}} \int_{\Lambda} g (\lambda, X_s) \, M (d \lambda, ds) \Big] 
		\\&= u_m (x) + E^P \Big[ \int_0^{t \wedge \tau_{D_m}} H_m (X_s, \nabla u_m (X_s), \nabla^2 u_m (X_s)) \, ds \Big]
		\\&\leq u_m (x) + E^P \Big[ \int_0^{\tau_{D_m}} H_m (X_s, \nabla u_m (X_s), \nabla^2 u_m (X_s))_+ \, ds \Big].
	\end{split}
	\end{equation}
	From Lemma~\ref{lem: main approx}~(3) and the Krylov estimate~\eqref{eq: krylov estimate}, we obtain
	\begin{equation} \label{eq: Krylov estimate lem 4.4}
		\begin{split}
	E^P \Big[ \int_0^{\tau_{D_m}} H_m (X_s&, \nabla u_m (X_s), \nabla^2 u_m (X_s))_+ \, ds \Big] 
	\\&\leq N \| H_m (\, \cdot \, , \nabla u_m, \nabla^2 u_m)_+ \|_{L_p (D_m)} = 0.
	\end{split}
	\end{equation} 
	Using once again the Krylov estimate~\eqref{eq: krylov estimate} and \(\sup_{\lambda \in \Lambda} |g (\lambda, x)| \leq \mathfrak{g} (x)\) by Standing Assumption~\ref{SA: main2}, we notice that
	\[
	E^P \Big[ \int_0^{\tau_{D_m}} \int_\Lambda |g (\lambda, X_s)| \, M (d \lambda, ds) \Big] \leq N \| \mathfrak{g} \|_{L_p (D)}.
	\] 
	Thus, as \(u_m\) is bounded on \(\oD_m\), the dominated convergence theorem yields that 
	\begin{align*}
		E^P \Big[ &u_m (X_{\tau_{D_m}}) + \int_0^{\tau_{D_m}} \int_\Lambda g (\lambda, X_s) \, M(d \lambda, ds) \Big] 
		\\&= \lim_{t \to \infty} E^P \Big[ u_m (X_{t \wedge \tau_{D_m}}) + \int_0^{t \wedge \tau_{D_m}} \int_\Lambda g (\lambda, X_s) \, M(d \lambda, ds) \Big] 
		\leq u_m (x), 
	\end{align*}
	where we used \eqref{eq: main inequality main step} and \eqref{eq: Krylov estimate lem 4.4} for the inequality. The proof is complete.
\end{proof}

\subsection{Extracting a convergent subsequence} \label{sec: subsequence}
Fix \(x \in D\) for the remainder of this proof. Let \(\J \in \mathbb{N}\) be such that \(x \in D_m\) for all \(m \geq \J\). 
Further, let \((P_m)_{m \geq \J}\) be the sequence from Lemma~\ref{lem: stochastic lower bound}.
By Lemma~\ref{lem: main approx}~(4), it satisfies
\begin{align} \label{eq: lower bound approximation stochastic}
u (x) =  \lim_{m \to \infty} u_m (x) \geq \limsup_{m \to \infty} E^{P_m} \Big[ u_m (X_{\tau_{D_m}}) + \int_0^{\tau_{D_m}} \int_\Lambda g (\lambda, X_s) \,  M(d \lambda, ds) \Big]. 
\end{align} 

\begin{lemma}\label{lem: tight + closed}
	The sequence \((P_m)_{m \geq \J}\) is tight and any of its accumulation points is in \(\cR (x)\). 
\end{lemma}

\begin{proof}
	Tightness of the marginal distributions \((P_m \circ X^{-1})_{m \geq \J}\) follows from \cite[Corollary~1.2]{krylov_21} and Kolmogorov's tightness criterion \cite[Theorem~23.7]{Kallenberg}. As the space \(\mathbb{M}\) is compact, we conclude that \((P_m)_{m \geq \J}\) is tight.
	We need to prove that any accumulation point is in \(\cR (x)\). To simplify the notation, assume that \(P_m \to P\). We have to establish the initial value condition and the martingale properties from the definition of \(\cK (x)\). For the initial value condition, notice that \(P_m (X_0 = x) = 1\) implies \(P (X_0 = x) = 1\), by the continuity of \(Q \mapsto Q \circ X^{-1}_0\). We now discuss the martingale properties. 
	We prepare the main proof with some technical observations.
	
	{\em Step 1:} 
		As the coefficients \((a_m, b_m)\) associated with \(P_m\) satisfy the properties from Standing Assumption~\ref{SA: main1}, the Krylov estimate \cite[Theorem~2.7]{krylov_AOP_21} yields, for every Borel function \(\ell \colon \bR^d \to \bR\) and \(T > 0\), that
	\begin{align} \label{eq: krylov Pn}
	E^{P_m} \Big[ \int_0^{T} |\ell (X_r)| \, dr \Big] \leq N \sqrt{T} \, \| \ell \|_{L_d (\bR^d)}, 
	\end{align} 
	where \(N = N (d, \delta, \| \b \|_{L_d (\bR^d)}) > 0\).
	For non-negative lower semicontinuous \(\ell\) this estimate transfers immediately to the limit \(P\) by \cite[Theorem~A.3.12]{DE_97}. We now explain that it also holds for \(P\) and any Borel function \(\ell\), i.e., 
	\begin{align} \label{eq: krylov convergence estimate}
	E^P \Big[ \int_0^{T} | \ell (X_r) | \, dr \Big] \leq N \sqrt{T} \, \| \ell \|_{L_d (\bR^d)}. 
	\end{align} 
    To see this, let \(\ell \in L_d (\bR^d)\) be a Borel function and \(\varepsilon > 0\). Then, the Vitali--Carath\'eodory theorem \cite[p. 56]{rudin_87} yields the existence of a non-negative lower semicontinuous function \(\ell_l\) such that \(|\ell |^d \leq \ell_l\) and 
    \[
    \int_{\bR^d} (\ell_l (z) - | \ell (z) |^d ) \, dz < \varepsilon^d. 
    \]
    Then, as \eqref{eq: krylov convergence estimate} holds for the lower semicontinuous function \(\ell_l^{1/d}\), we obtain that 
    \begin{align*}
        E^P \Big[ \int_0^{T} | \ell (X_r) | \, dr \Big] &\leq E^P \Big[ \int_0^{T}  \ell_l^{1/d} (X_r) \, dr \Big] 
        \\&\leq N \sqrt{T} \, \| \ell_l^{1/d} \|_{L_d (\bR^d)}
        \leq N \sqrt{T} \, \| \ell \|_{L_d (\bR^d)} + N \sqrt{T} \, \varepsilon. \phantom \int 
    \end{align*}
    Since \(\varepsilon > 0\) was arbitrary, we conclude that \eqref{eq: krylov convergence estimate} holds for every Borel function \(\ell\). 
	
	{\em Step 2:}
	Fix an arbitrary \(\varphi \in C^2_c (\bR^d)\) and let \(C (\varphi)\) be the process from \eqref{eq: test process}, i.e., 
	\[
		C (\varphi) = \varphi (X) - \varphi (X_0) - \int_{0}^{\cdot} \int_{\Lambda}  L (\lambda, X_s, \nabla \varphi (X_s), \nabla^2 \varphi (X_s)) \, M (d \lambda, ds)
	\]
	with
	\begin{align*}
		L (\lambda, x, \nabla \varphi (x), \nabla^2 \varphi (x)) &=  \langle b (\lambda, x), \nabla \varphi (x) \rangle + \tfrac{1}{2} \on{tr} \big[ a (\lambda, x) \nabla^2 \varphi (x) \big].
	\end{align*}
	By a standard monotone class argument, the \((\mathcal{G}_t)_{t \geq 0}\)-\(P\)-martingale property of \(C (\varphi)\) can be rephrased as follows: for all \(0 \leq s < t < \infty\) and all \(\mathcal{G}_s\)-measurable \(\psi \in C_b (\Omega \times \mathbb{M})\), it holds that
	\begin{align} \label{eq: supermartingale to show}
		E^P \Big[ \psi (\, \cdot \,) \big( C_{t} (\varphi) - C_{s} (\varphi) \big) \Big] = 0.
	\end{align} 
As the test function \(\varphi\) is compactly supported, there exists a compact set \(K \subset \bR^d\) such that \[
L (\lambda, x, \nabla \varphi (x), \nabla^2 \varphi (x)) = 0\] for all \(x \not \in K\). 
Using this observation, we may argue as in the proof for \eqref{eq: approx in Lp} and obtain that, for every \(\varepsilon > 0\), there exists a bounded continuous function \(\ell_\varepsilon \colon \Lambda \times \bR^d \to \bR\) such that 
\begin{align} \label{eq: g eps}
\Big\| \sup_{\lambda \in \Lambda} | L (\lambda, \cdot \, , \nabla \varphi, \nabla^2 \varphi) - \ell_\varepsilon (\lambda, \, \cdot \,) | \Big\|_{L_d (\bR^d)} \leq \varepsilon. 
\end{align}
	As \(P_m \in \cK_m (x)\), we have 
	\begin{align} \label{eq: leq 0}
E^{P_m} \Big[ \psi (\, \cdot \, )\Big( \varphi (X_{t}) - \varphi (X_{s}) - \int_{s}^{t} \int_\Lambda L_m (\lambda, X_r, \nabla \varphi (X_r), \nabla^2 \varphi (X_r)) \,M(d \lambda, dr) \Big) \Big] = 0, 
	\end{align} 
	where
		\begin{align*}
		L_m (\lambda, x, \nabla \varphi (x), \nabla^2 \varphi (x)) &:=  \langle b_m (\lambda, x), \nabla \varphi (x) \rangle + \tfrac{1}{2} \on{tr} \big[ a_m (\lambda, x) \nabla^2 \varphi (x) \big].
	\end{align*}
	Moreover, by the boundedness and continuity of \(\ell_\varepsilon\) and \(\psi\), the continuous mapping theorem yields that 
	\begin{align} \label{eq: convergence with continuous input}
		\lim_{m \to \infty} E^{P_m} \Big[ \psi \int_{s}^{t} \int_\Lambda  \ell_\varepsilon (\lambda, X_r) \, M (d \lambda, dr) \Big] = E^P \Big[ \psi \int_{s}^{t} \int_\Lambda \ell_\varepsilon (\lambda, X_r) \, M (d \lambda, dr) \Big]. 
	\end{align}
Using the Krylov estimates \eqref{eq: krylov Pn} and \eqref{eq: krylov convergence estimate}, and \eqref{eq: g eps}, we find that  
	\begin{equation} \label{eq: long estimate lem 4.4}
		\begin{split}
		\Big| E^P &\Big[ \psi \int_{s}^{t} \int_\Lambda L (\lambda, X_r, \nabla \varphi(X_r), \nabla^2 \varphi (X_r)) \, M (d\lambda, dr) \Big] \\&\hspace{3cm}- E^{P_m} \Big[ \psi \int_{s}^{t} \int_\Lambda L_m (\lambda, X_r, \nabla \varphi(X_r), \nabla^2 \varphi (X_r)) \, M (d \lambda, dr) \Big] \Big|
		\\&\leq  
		\Big| E^{P} \Big[ \psi \int_{s}^{t} \int_\Lambda (L (\lambda, X_r, \nabla \varphi(X_r), \nabla^2 \varphi (X_r)) - \ell_\varepsilon (\lambda, X_r)) \, M (d \lambda, dr) \Big] \Big| 
		\\&\qquad \qquad + \Big| E^{P}  \Big[ \psi \int_{s}^{t} \int_\Lambda  \ell_\varepsilon (\lambda, X_r) \, M (d \lambda, dr) \Big] - E^{P_m} \Big[ \psi \int_{s}^{t} \int_\Lambda \ell_\varepsilon (\lambda, X_r) \, M (d \lambda, dr) \Big] \Big| 
		\\& \qquad \qquad + \Big| E^{P_m} \Big[ \psi \int_{s}^{t} \int_\Lambda (\ell_\varepsilon (\lambda, X_r) - L_m (\lambda, X_r, \nabla \varphi(X_r), \nabla^2 \varphi (X_r))) \, M (d \lambda, dr) \Big] \Big| 
		\\&\leq N \sqrt{t} \, \|\psi\|_\infty \Big\| \sup_{\lambda \in \Lambda} | L (\lambda, \cdot \, , \nabla \varphi, \nabla^2 \varphi) - \ell_\varepsilon (\lambda, \, \cdot \,) | \Big\|_{L_d (\bR^d)} 
		\\&\qquad \qquad + \Big| E^{P}  \Big[ \psi \int_{s}^{t} \int_\Lambda  \ell_\varepsilon (\lambda, X_r) \, M (d \lambda, dr) \Big] - E^{P_m} \Big[ \psi \int_{s}^{t} \int_\Lambda \ell_\varepsilon (\lambda, X_r) \, M (d \lambda, dr) \Big] \Big| 
		\\&\qquad \qquad + N \sqrt{t} \, \|\psi\|_\infty \Big\| \sup_{\lambda \in \Lambda} | L (\lambda, \cdot \, , \nabla \varphi, \nabla^2 \varphi) - \ell_\varepsilon (\lambda, \, \cdot \,) | \Big\|_{L_d (\bR^d)}  \phantom \int 
        \\&\qquad \qquad + \textup{const. } N \sqrt{t}\, \|\psi\|_\infty \| (1 + \b ) \1_{G_m} \|_{L_d (\bR^d)} \phantom \int 
        \\&\leq  2 N \sqrt{t} \, \| \psi \|_\infty \, \varepsilon + \Big| E^{P}  \Big[ \psi \int_{s}^{t} \int_\Lambda  \ell_\varepsilon (\lambda, X_r) \, M (d \lambda, dr) \Big] - E^{P_m} \Big[ \psi \int_{s}^{t} \int_\Lambda \ell_\varepsilon (\lambda, X_r) \, M (d \lambda, dr) \Big] \Big| 
        \\&\qquad \qquad + \textup{const. } N \sqrt{t}\, \|\psi\|_\infty \| (1 + \b ) \1_{G_m} \|_{L_d (\bR^d)}, \phantom \int 
	\end{split}
	\end{equation} 
	where the non-specified constant only depends on \(d, \delta\) and \(\varphi\).
	Lemma~\ref{lem: main approx}~(2) and the dominated convergence theorem (for convergence in measure) yield that 
	\[
	\| (1 + \b ) \1_{G_m} \|_{L_d (\bR^d)} \to 0, \quad m \to \infty. 
	\] 
	Using this and \eqref{eq: convergence with continuous input}, letting first \(m \to \infty\) and then \(\varepsilon \to 0\), we obtain that 
	\begin{equation} \label{eq: second term convergence} \begin{split}
	E^P \Big[ \psi \int_{s}^{t} &\int_\Lambda L (\lambda, X_r, \nabla \varphi(X_r), \nabla^2 \varphi (X_r)) \, M (d\lambda, dr) \Big]
	\\&= \lim_{m \to \infty} E^{P_m} \Big[ \psi \int_{s}^{t} \int_\Lambda L_m (\lambda, X_r, \nabla \varphi(X_r), \nabla^2 \varphi (X_r)) \, M (d \lambda, dr) \Big].
	\end{split}
	\end{equation} 
	Finally, by the boundedness and continuity of \(\psi\) and \(\varphi\), the continuous mapping theorem yields that 
	\begin{align} \label{eq: first term convergence}
	\lim_{m \to \infty} E^{P_m} \Big[ \psi (\, \cdot \,) \Big( \varphi (X_{t}) - \varphi(X_{s}) \Big) \Big] = E^{P} \Big[ \psi (\, \cdot \,) \Big( \varphi (X_{t}) - \varphi (X_{s}) \Big) \Big].
	\end{align} 
	To summarize, \eqref{eq: leq 0}, \eqref{eq: second term convergence} and \eqref{eq: first term convergence} yield \eqref{eq: supermartingale to show}, which completes the proof. 
\end{proof}

\begin{lemma} \label{lem: continuity} 
	Let \((P_{m_n})_{n = 1}^\infty\) be a convergent subsequence of \((P_m)_{m \geq \J}\) and denote its limit by \(P\). Then, 
	\begin{align*}
	\lim_{n \to \infty} E^{P_{m_n}} \Big[ u_{m_n} (X_{\tau_{D_{m_n}}}&) + \int_0^{\tau_{D_{m_n}}}  \int_\Lambda g (\lambda, X_s) \, M (d \lambda, ds) \Big] 
	\\&= E^P \Big[ u (X_{\tau_D}) + \int_0^{\tau_D} \int_\Lambda g (\lambda, X_s) \, M (d\lambda, ds) \Big]. 
	\end{align*}
\end{lemma} 

\begin{proof}
	Set \(F := \{ \omega \in \Omega \colon \tau_D (\omega) = \tau_{\overline{D}} (\omega) < \infty \}\), where \(\tau_{\oD} (\omega) := \inf \{t \geq 0 \colon \omega (t) \not \in \oD \}\). By Lemma \ref{lem: tight + closed} and \cite[Lemma~3.5]{C_25_SIMA}, \(P (F) = 1\) and \(\omega \mapsto \tau_D (\omega)\) is continuous at every point in \(F\). Let \((\omega_n)_{n = 1}^\infty \subset \Omega\) be a sequence with \(\omega_n (0) = x\) and \(\omega_n \to \omega \in F\). As \(D_{m_n} \subset D\),
	\[
	\limsup_{n \to \infty} \tau_{D_{m_n}} (\omega_n) \leq \limsup_{n \to \infty} \tau_D (\omega_n) = \tau_D (\omega) < \infty.
	\]
	Conversely, take \(t < \tau_D (\omega)\). Then, \(I_t := \{ \omega (s) \colon s \in [0, t] \}\) is a compact subset of \(D = \bigcup_{K = 1}^\infty D_K\), and hence \(I_t \subset D_K\) for some \(K \in \mathbb{N}\). Let \(\varepsilon := \on{dist} (I_t, \bR^d \setminus D_K) > 0\). Since \(\omega_n \to \omega\) locally uniformly and \(m_n \to \infty\), there exists an \(n_0 \in \mathbb{N}\) such that, for all \(n \geq n_0\), we have \(m_n \geq K\) and
	\[
	\sup_{s \in [0, t]} \| \omega_n (s) - \omega(s)\| < \varepsilon,
	\]
	which implies \(\{ \omega_n (s) \colon s \in [0, t] \} \subset D_K \subset D_{m_n}\), and hence \(\tau_{D_{m_n}} (\omega_n) \geq t\). Letting \(t \nearrow \tau_D (\omega)\) shows that \(\tau_{D_{m_n}} (\omega_n) \to \tau_D (\omega)\). As \(\omega_n \to \omega\) locally uniformly, this yields \(y_n := \omega_n (\tau_{D_{m_n}} (\omega_n)) \to y := \omega (\tau_D (\omega))\). Since \(y_n \in \oD_{m_n}\) for all sufficiently large \(n\), Lemma~\ref{lem: main approx}~(4) and the continuity of \(u\) on \(\oD\) yield
	\[
	| u_{m_n} (y_n) - u (y) | \leq \sup_{z \in \oD_{m_n}} | u_{m_n} (z) - u (z) | + | u (y_n) - u (y) | \to 0, \quad n \to \infty.
	\]
	Now, the uniform boundedness 
	\[
	\sup_{m \in \mathbb{N}} \sup_{z \in \oD_m} | u_m (z) | < \infty, 
	\] 
	which is implied by Lemma~\ref{lem: main approx}~(4), and the extended continuous mapping theorem \cite[Theorem~5.27]{Kallenberg} yield that 
	\[
	E^{P_{m_n}} \big[ u_{m_n} (X_{\tau_{D_{m_n}}}) \big] \to E^P \big[ u (X_{\tau_D}) \big]. 
	\] 
	For the second term, take \(\xi > 0\). By \eqref{eq: approx in Lp}, there exists a function \(g_\xi \in C_b (\Lambda \times \bR^d)\) such that 
	\[
	\Big\| \sup_{\lambda \in \Lambda} | g (\lambda, \, \cdot \,) - g_\xi (\lambda, \, \cdot \,)| \Big \|_{L_p (D)} \leq \xi.
	\] 	
	Using the Krylov estimate~\eqref{eq: krylov estimate}, we find that 
	\begin{align*}
		\Big| E^{P_{m_n}} &\Big[  \int_0^{\tau_{D_{m_n}}} \int_\Lambda g (\lambda, X_s) \, M(d \lambda, ds) \Big] - E^P \Big[ \int_0^{\tau_D} \int_\Lambda g (\lambda, X_s) \, M(d \lambda, ds) \Big] \Big| 
		\\&\leq 2 N \Big\| \sup_{\lambda \in \Lambda} | g (\lambda, \, \cdot \,) - g_\xi (\lambda, \, \cdot \,)| \Big \|_{L_p (D)}
		\\&\qquad + \Big| E^{P_{m_n}} \Big[ \int_0^{\tau_{D_{m_n}}} \int_\Lambda g_\xi (\lambda, X_s) \, M (d \lambda, ds) \Big] - E^P \Big[ \int_0^{\tau_D} \int_\Lambda g_\xi (\lambda, X_s) \, M (d \lambda, ds) \Big] \Big|
			\\&\leq 2 N \xi + \Big| E^{P_{m_n}} \Big[ \int_0^{\tau_{D_{m_n}}} \int_\Lambda g_\xi (\lambda, X_s) \, M (d \lambda, ds) \Big] - E^P \Big[ \int_0^{\tau_D} \int_\Lambda g_\xi (\lambda, X_s) \, M (d \lambda, ds) \Big] \Big|.
		\end{align*}
			For \((\omega_n, \nu_n) \to (\omega, \nu) \in F \times \mathbb{M}\), we observe that 
		\begin{align*}
			\Big| &	\int_0^{\tau_{D_{m_n}} (\omega_n)}  \int_\Lambda g_\xi (\lambda, \omega_n (s)) \, \nu_n (d\lambda, ds) - \int_0^{\tau_D (\omega)} \int_\Lambda g_\xi (\lambda, \omega (s)) \, \nu (d \lambda, ds) \Big| 
			\\&\quad\leq | \tau_{D_{m_n}} (\omega_n) - \tau_D (\omega) |  \, \| g_\xi \|_\infty + \tau_D (\omega) \sup_{\substack{\lambda \in \Lambda \\ s \in [0, \tau_D (\omega)]}} | g_\xi (\lambda, \omega_n (s)) - g_\xi (\lambda, \omega (s)) | 
			\\&\quad \qquad+ \Big|	\int_0^{\tau_{D} (\omega)} \int_\Lambda g_\xi (\lambda, \omega (s)) \, (\nu_n - \nu) (d\lambda, ds) \Big| \to 0, 
		\end{align*} 
		as \(n \to \infty\). Here, the first term vanishes as \(\tau_{D_{m_n}} (\omega_n) \to \tau_D (\omega)\), as shown above, and the second term vanishes by the uniform continuity of \(g_\xi\) on compact subsets of \(\Lambda \times \bR^d\). For the last term, let \(\eta > 0\) and take a continuous function \(\chi_\eta \colon \bR_+ \to [0, 1]\) with \(\chi_\eta = 1\) on \([0, \tau_D (\omega)]\) and \(\chi_\eta = 0\) on \([\tau_D (\omega) + \eta, \infty)\). As \((s, \lambda) \mapsto \chi_\eta (s) g_\xi (\lambda, \omega (s))\) is continuous with compact support in \(\bR_+ \times \Lambda\), the vague convergence \(\nu_n \to \nu\) yields
		\[
		\int_0^\infty \int_\Lambda \chi_\eta (s) g_\xi (\lambda, \omega (s)) \, (\nu_n - \nu) (d \lambda, ds) \to 0, \quad n \to \infty.
		\]
		Furthermore, as all elements of \(\m\) have the Lebesgue measure as projection to \(\bR_+\),
		\[
		\Big| \int_0^\infty \int_\Lambda \big( \1_{[0, \tau_D (\omega)]} (s) - \chi_\eta (s) \big) g_\xi (\lambda, \omega (s)) \, \mu (d \lambda, ds) \Big| \leq \eta \, \| g_\xi \|_\infty, \quad \mu \in \{\nu, \nu_1, \nu_2, \dots\}.
		\]
		Letting first \(n \to \infty\) and then \(\eta \searrow 0\) shows that the last term vanishes as well.
		Hence, the extended continuous mapping theorem \cite[Theorem~5.27]{Kallenberg} yields that, weakly as \(n \to \infty\), 
		\[
		P_{m_n} \circ \Big( \int_0^{\tau_{D_{m_n}}} \int_\Lambda g_\xi (\lambda, X_s) \, M (d \lambda, ds) \Big)^{-1} \to P \circ \Big( \int_0^{\tau_D} \int_\Lambda g_\xi (\lambda, X_s) \, M (d \lambda, ds) \Big)^{-1}.
		\] 
		 We deduce from \cite[Corollary~2.7]{krylov_21} that 
		 \[
		 \sup_{n \in \mathbb{N}} E^{P_{m_n}} \big[ \tau_D^2 \big] < \infty. 
		 \] 
		 Thus, because 
		\[
	\Big|	\int_0^{\tau_{D_{m_n}}} \int_\Lambda g_\xi (\lambda, X_s) \, M (d \lambda, ds) \Big| \leq \tau_{D} \, \sup_{\substack {\lambda \in \Lambda \\ y \in \oD}} | g_\xi (\lambda, y) |, 
		\] 
		the sequence 
		\[
		\Big( P_{m_n} \circ \Big( \int_0^{\tau_{D_{m_n}}} \int_\Lambda g_\xi (\lambda, X_s) \, M (d \lambda, ds) \Big)^{-1} \Big)_{n = 1}^\infty
		\] 
is uniformly integrable and we obtain that 
		\[
		\Big| E^{P_{m_n}} \Big[ \int_0^{\tau_{D_{m_n}}} \int_\Lambda g_\xi (\lambda, X_s) \, M (d\lambda, ds) \Big] - E^P \Big[ \int_0^{\tau_D} \int_\Lambda g_\xi (\lambda, X_s) \, M (d \lambda, ds) \Big] \Big| \to 0, \quad n \to \infty.
		\]
		Summarizing, taking first \(n \to \infty\) and afterwards \(\xi \to 0\), we conclude that 
		\[
		E^{P_{m_n}} \Big[ \int_0^{\tau_{D_{m_n}}} \int_\Lambda g (\lambda, X_s) \, M (d \lambda, ds) \Big] \to E^P \Big[ \int_0^{\tau_D} \int_\Lambda g (\lambda, X_s) \, M (d\lambda, ds) \Big], 
		\] 
		completing the proof. 
\end{proof}

\subsection{Conclusion}
For \(x \in \partial D\), we have \(v (x) = f (x) \leq u (x)\) by the definition of an \(L_p\)-viscosity supersolution to \eqref{eq: PDE}. For \(x \in D\), Lemma~\ref{lem: tight + closed} yields a subsequence of \((P_m)_{m \geq \J}\) that converges weakly to some \(P \in \cR (x)\). As the \(\limsup\) in \eqref{eq: lower bound approximation stochastic} dominates the limit along this subsequence, Lemma~\ref{lem: continuity} and \(u \geq f\) on \(\partial D\) yield that
\begin{equation} \label{eq: final conclusion}
	\begin{split} 
	u (x) &\geq \limsup_{m \to \infty} E^{P_m} \Big[ u_m (X_{\tau_{D_m}}) + \int_0^{\tau_{D_m}} \int_\Lambda g (\lambda, X_s) \, M (d \lambda, ds) \Big] \\&\geq E^{P} \Big[ f (X_{\tau_D}) + \int_0^{\tau_D} \int_\Lambda g (\lambda, X_s) \, M (d \lambda, ds) \Big] 
	\\&\geq \inf_{Q \in \cR (x)} E^Q \Big[ f (X_{\tau_D}) + \int_0^{\tau_D}\int_\Lambda g (\lambda, X_s) \, M (d\lambda, ds) \Big] = v (x).
\end{split}
\end{equation} 
As we already know from Section~\ref{sec: VP} that \(v\) is an \(L_p\)-viscosity solution to \eqref{eq: PDE}, and because \(u\) was an arbitrary \(L_p\)-viscosity supersolution (see Standing Assumption~\ref{SA: p, u}), the proof of Theorem~\ref{theo: main1} is complete. \qed

	\section*{Acknowledgment}
	
The authors used ChatGPT~6 (OpenAI) and Claude Opus~5.5 (Anthropic) to identify typographical and minor mathematical errors, to assist with the literature review, and to discuss and draft formalizations of some of the authors' own proof ideas for Proposition~\ref{prop: representation}, Lemma~\ref{lem: main approx}, and Appendix~\ref{ex: non unique}. All suggestions were carefully reviewed, revised where necessary, and verified by the authors, who take full responsibility for the content and mathematical correctness of the paper.

	\appendix 
	
	\section{Example of a semilinear Dirichlet problem \\with infinitely many viscosity solutions} \label{ex: non unique}
	
		Providing some intuition for the construction, we modify Nadirashvili's linear example from \cite{zbMATH01140123} by a small controlled drift, producing a nonlinear first order term. In case the influence of the drift is small enough, the nonlinear equation is close enough to the linear equation so that its non-uniqueness propagates. Let us provide the details. 
		
		We work in the framework from Example~\ref{ex: coefficients}. More specifically, let \(D = B_1 (0)\) be the unit ball in \(\bR^d\) with \(d \geq 3\) and take \(A \colon D \to \mathbb{S}_{\delta^*}^d\) and \(f \in C^2 (\partial D)\) as in the {\em linear} example from \cite[Theorem on p. 538]{zbMATH01140123}. Here, \(\delta^* \in (0, 1)\) is some constant that comes from the coefficient \(A\) in \cite{zbMATH01140123}. We extend the coefficient \(A\) from \(D\) to \(\bR^d\) by setting \(A (x) := \bar{A}\) for \(x \in \mathbb{R}^d \setminus D\) and an arbitrary matrix \(\bar{A} \in \mathbb{S}_{\delta^*}^d\). Let \(\xi \in C^\infty_c (\bR^d)\) with \(\xi = 1\) on \(D\). Then, similarly to Example~\ref{ex: coefficients}, for \(\Lambda := [-1, 1]\), \(\alpha \in (0, 1)\), \(e \in \partial B_1 (0)\) and \(\varepsilon > 0\), we set 
		\[
		a (\lambda, x) := (1 - \alpha \lambda) A (x), \quad b (\lambda, x) := \lambda \varepsilon \xi (x) e, 
		\]
		where \(a (\lambda, x) \in \mathbb{S}_\delta^d\) with \(\delta := (1 - \alpha) \delta^*\). Throughout this example, we take \(\varepsilon \in (0, 1]\). Then, \(\| b (\lambda, x) \| \leq |\xi (x)|\) for all \((\lambda, x) \in \Lambda \times \bR^d\), i.e., Standing Assumption~\ref{SA: main1}~(c) holds with \(\b := |\xi|\) independently of \(\varepsilon\). In particular, the constant \(d_0\) from Theorem~\ref{theo: main1} can be chosen independently of \(\varepsilon \in (0, 1]\).
		In this case, straightforward considerations show that the Dirichlet problem \eqref{eq: PDE} with \(g \equiv 0\) has the same \(L_p\)-viscosity sub- and supersolutions as
		\begin{align} \label{eq: PDE non unique}
			\begin{cases} \displaystyle
				- \tfrac{1}{2} \on{tr} \big[ A (x) \nabla^2 u (x) \big] + \frac{\varepsilon}{1 - \alpha^2} \Big( | \langle e, \nabla u (x) \rangle | - \alpha \langle e, \nabla u (x) \rangle \Big) = 0, & x \in D, \\ u (x) = f (x), & x \in \partial D.
			\end{cases} 
		\end{align} 
		For sufficiently small \(\varepsilon > 0\), we now show that the nonlinear Dirichlet problem \eqref{eq: PDE} has infinitely many \(L_p\)-viscosity solutions for all \(p \in (d_0, d)\).
		
		By \cite[Theorem on p. 538]{zbMATH01140123}, after a cutoff and diagonal approximation, if necessary, there exist two sequences \((A^k_n)_{n = 1}^\infty \colon \bR^d \to \mathbb{S}_{\delta^*}^d \subset \mathbb{S}_\delta^d\), \(k = 1, 2\), of smooth functions and classical solutions \((v^1_n)_{n = 1}^\infty\) and \((v^2_n)_{n = 1}^\infty\) to the linear Dirichlet problem 
		\[
		\begin{cases} \displaystyle \tfrac{1}{2} \on{tr} \big[ A^k_n (x) \nabla^2 v^k_n (x)\big] = 0, & x \in D, \\ v^k_n (x) = f (x), & x \in \partial D, \end{cases}
		\]
		such that \(A^k_n \to A\) a.e., \(v^k_n \to v^k\) in \(C (\oD)\) and \(v^1 (x_0) \not = v^2 (x_0)\) for some \(x_0 \in D\).
		By \cite[Theorem~1.2, Example~1.2]{krylov_CPDE_20}, there are Sobolev solutions \(w^k_n \in W^2_{p, \on{loc}} (D) \cap C (\oD)\) to the nonlinear Dirichlet problems 
		\begin{align} \label{eq: non-linear approx}
			\begin{cases} \displaystyle
				- \tfrac{1}{2} \on{tr} \big[ A^k_n (x) \nabla^2 w^k_n (x) \big] + \frac{\varepsilon}{1 - \alpha^2} \Big( | \langle e, \nabla w^k_n (x) \rangle | - \alpha \langle e, \nabla w^k_n (x) \rangle \Big) = 0, & x \in D, \\ w^k_n (x) = f (x), & x \in \partial D.
			\end{cases} 
		\end{align} 
		As the coefficient \(a^k_n (\lambda, x) := (1 - \alpha \lambda) A^k_n (x)\) is uniformly continuous with values in \(\mathbb{S}_{\delta}^d\), \cite[Proposition~4.5]{KS_22} and \cite[Proposition~2.7]{C_25_SIMA} yield that \(w^k_n\) has a stochastic control representation as a value function with cost \(f (X_{\tau_D})\).
		This observation allows us to apply the stability result \cite[Proposition~2.9]{C_25_SIMA}, which shows that, along a subsequence, \((w^k_n)_{n = 1}^\infty\) converges in \(C (\oD)\) and its limit point \(w^k\) is an \(L_p\)-viscosity solution to the Dirichlet problem \eqref{eq: PDE}. We now show that \(w^1 (x_0) \not = w^2 (x_0)\) for \(\varepsilon > 0\) small enough. 
		
		To prove this, we show that
		\begin{align} \label{eq: to show example}
			\sup_{x \in \oD} | w^k - v^k | \leq \frac{1}{3} | v^1 (x_0) - v^2 (x_0) |, \quad k = 1, 2, 
		\end{align} 
		once \(\varepsilon > 0\) is taken small enough.
		This entails that \(w^1 (x_0) \not = w^2 (x_0)\), because 
		\begin{align*} 
			| w^1 (x_0) - w^2 (x_0) | &\geq | v^1 (x_0) - v^2 (x_0) | - | w^1 (x_0) - v^1 (x_0) | - | v^2 (x_0) - w^2 (x_0)| 
			\\&\geq \frac{1}{3} | v^1 (x_0) - v^2 (x_0) | > 0.  
		\end{align*} 
		In summary, it remains to establish \eqref{eq: to show example}. 
		
		For a fixed \(\bar{x} \in D = B_1 (0)\), let \(Y\) be a solution process to the stochastic differential equation 
		\[
		dY_t = \sqrt{A^k_n (Y_t)} \, dW_t, \quad Y_0 = \bar{x}, 
		\]
		where \(W\) is a \(d\)-dimensional Brownian motion.
		As \(A^k_n\) is smooth with values in \(\mathbb{S}_\delta^d\), strong existence for this equation is guaranteed by standard existence theory (see, e.g., Chapter 5 in \cite{KaraShre}). Let \(\tau_h := \inf \{t \geq 0 \colon \|Y_t\| \geq 1 - 1/h\} \wedge h\).
		Using Krylov's It\^o formula \cite[Theorem~1.3]{krylov_AOP_21} and \eqref{eq: non-linear approx}, we obtain that 
		\begin{align*}
			v^k_n (\bar{x}) - w^k_n (\bar{x}) = \E \big[ v^k_n (Y_{\tau_h}) - w^k_n (Y_{\tau_h}) \big] + \frac{\varepsilon}{1 - \alpha^2} \, \E \Big[ \int_0^{\tau_h} \Big( | \langle e, \nabla w^k_n (Y_t) \rangle | - \alpha \langle e, \nabla w^k_n (Y_t) \rangle \Big) \, dt \Big], 
		\end{align*}
		which implies that 
		\begin{align} \label{eq: main to control}
			| v^k_n (\bar{x}) - w^k_n (\bar{x}) | \leq \Big| \E \big[ v^k_n (Y_{\tau_h}) - w^k_n (Y_{\tau_h}) \big] \Big| + \frac{2 \varepsilon}{1 - \alpha^2} \E \Big[ \int_0^{\tau_h} \| \nabla w^k_n (Y_t) \| \, dt \Big]. 
		\end{align}
		As \(h \to \infty\), we have 
		\[
		\E \big[ v^k_n (Y_{\tau_h}) - w^k_n (Y_{\tau_h}) \big] \to \E \big[ v^k_n (Y_{\tau_D (Y)}) - w^k_n (Y_{\tau_D (Y)}) \big] = \E \big[ f (Y_{\tau_D (Y)}) - f (Y_{\tau_D (Y)}) \big] = 0.
		\]
		Hence, in case we control the term 
		\[
		\E \Big[ \int_0^{\tau_h} \| \nabla w^k_n (Y_t) \| \, dt \Big]
		\]
		uniformly in \(n\) and \(h\), we may choose \(\varepsilon > 0\) in \eqref{eq: main to control} small enough such that \eqref{eq: to show example} holds. This is the task for the remainder of this proof. 
		
		Let \(m_f := \min_{\partial D} f\) and \(M_f := \max_{\partial D} f\). Recalling that \(w^k_n\) has a stochastic control representation as a value function with cost \(f (X_{\tau_D})\), we get that  
		\begin{align}  \label{eq: upper lower bound} 
			m_f \leq w^k_n \leq M_f \quad \text{on } \oD.
		\end{align}
		Using \eqref{eq: non-linear approx}, \eqref{eq: upper lower bound} and \(A^k_n (\, \cdot \,) \in \mathbb{S}_\delta^d\), a.e. on \(D\), we obtain that 
		\begin{equation} \label{eq: useful estimate} 
			\begin{split}
				\tfrac{1}{2} \on{tr} \big[ A^k_n \nabla^2 (w^k_n - m_f)^2 \big] &= \langle \nabla w^k_n, A^k_n \nabla w^k_n \rangle + \tfrac{2 \varepsilon}{1 - \alpha^2} (w^k_n - m_f) ( | \langle e, \nabla w^k_n \rangle | - \alpha \langle e, \nabla w^k_n \rangle )
				\\&\geq \langle \nabla w^k_n, A^k_n \nabla w^k_n \rangle \geq \delta \| \nabla w^k_n \|^2.
			\end{split}
		\end{equation} 
		By the Sobolev embedding theorem (see \cite[Theorem~4.12~Case C, p.~85]{adams_03}), we have 
		\[
		W^{2}_{p, \on{loc}} (D) \subset W^1_{q, \on{loc}} (D), \quad q = \frac{dp}{d - p} > 2p, \text{ since } p > d /2, 
		\] 
		which yields that \(\nabla^2 ((w^k_n - m_f)^2) = 2 \nabla w^k_n \otimes \nabla w^k_n + 2 (w^k_n - m_f) \nabla^2 w^k_n \in L_{p, \on{loc}} (D)\), using also the boundedness of \(w^k_n\) through its continuity.
		From this it is easily seen that \((w^k_n - m_f)^2 \in W^2_{p, \on{loc}} (D)\) and we may apply Krylov's It\^o formula \cite[Theorem~1.3]{krylov_AOP_21} to deduce from the estimate \eqref{eq: useful estimate} that
		\begin{equation} \label{eq: ex below bound} 
			\begin{split}
				\E \big[ ( w^k_n (Y_{\tau_h}) - m_f)^2 \big] - ( w^k_n (\bar{x}) - m_f)^2 &= \E \Big[ \int_0^{\tau_h} \tfrac{1}{2} \on{tr} \big[ A^k_n (Y_t) \nabla^2 ((w^k_n - m_f)^2) (Y_t) \big] \, dt \Big] 
				\\&\geq \delta \E \Big[ \int_0^{\tau_h} \| \nabla w^k_n (Y_t) \|^2 \, dt \Big].
			\end{split}
		\end{equation} 
		Combining \eqref{eq: upper lower bound} with \eqref{eq: ex below bound}, we find that 
		\[
		\E \Big[ \int_0^{\tau_h} \| \nabla w^k_n (Y_t) \|^2 \, dt \Big] \leq \frac{(M_f - m_f)^2}{\delta}.
		\]
		Returning to \eqref{eq: main to control}, letting \(h \to \infty\) and using H\"older's inequality, we conclude that 
		\[
		| v^k_n (\bar{x}) - w^k_n (\bar{x}) | \leq \frac{2 \varepsilon (M_f - m_f)}{\sqrt{\delta} (1 - \alpha^2)} \,  \sqrt{\E \big[ \tau_D \big]}. 
		\]
		By \cite[Theorem~2.8]{krylov_AOP_21}, there is a constant \(C (d, \delta^*) > 0\) such that \(\E [ \tau_D ] \leq C (d, \delta^*)\). Thus, we may take \(\varepsilon > 0\) small enough such that \eqref{eq: to show example} holds. 
		
		Summarizing, we proved that the Dirichlet problem \eqref{eq: PDE} has at least two different \(L_p\)-viscosity solutions. In fact, it has even infinitely many different \(L_p\)-viscosity solutions. To see this, let \(w_{\theta, n}\) be the Sobolev solution corresponding to the diffusion coefficient \(\theta A^1_n + (1 - \theta) A^2_n\) for \(\theta \in [0, 1]\), which exists by \cite[Theorem~1.2]{krylov_CPDE_20}. The stability and uniqueness results from \cite[Propositions~2.7, 2.9]{C_25_SIMA} show that \(\theta \mapsto w_{\theta, n} (x_0)\) is continuous. Let \((n (m))_{m = 1}^\infty\) be a subsequence such that
		\[
		\sup_{x \in \oD} \big| w^k_{n (m)} (x) - w^k (x) \big| \to 0, \quad m \to \infty, 
		\]
		for both \(k = 1, 2\). 
		For every \(r \in (w^1 (x_0) \wedge w^2 (x_0), w^1 (x_0) \vee w^2 (x_0))\), there exist an \(m_0 \in \mathbb{N}\) and \(\theta_{n (m)} \in [0, 1]\) with \(w_{\theta_{n (m)}, n (m)} (x_0) = r\) for all \(m \geq m_0\). 
		Notice that \(\theta A^1_n + (1 - \theta) A^2_n \in \mathbb{S}_{\delta^*}^d\) and that, for a.e. \(x \in D\), as \(m \to \infty\),  
		\[
		\| \theta_{n (m)} A^1_{n (m)} (x) + (1 - \theta_{n (m)}) A^2_{n (m)} (x) - A (x) \| \leq \| A^1_{n (m)} (x) - A (x) \| + \| A^2_{n (m)}  (x)- A (x) \| \to 0. 
		\]
		Hence, by \cite[Proposition~2.9]{C_25_SIMA}, the sequence \((w_{\theta_{n (m)}, n (m)})_{m = m_0}^\infty\) has an accumulation point \(w^r\) in \(C (\oD)\), which is an \(L_p\)-viscosity solution of \eqref{eq: PDE}. As \(w^r (x_0) = r\), the solutions \(w^r\), \(r \in (w^1 (x_0) \wedge w^2 (x_0), w^1 (x_0) \vee w^2 (x_0))\), are pairwise distinct.
	
	\section{Proof of Proposition~\ref{prop: representation}} \label{app: pf prop}
		Under the assumptions (i)--(iii), the set 
		\[ 
		\overline{\Theta} (x) := \Big\{ (A, q) \in \mathbb{S}^d \times \bR^d \colon \langle q, q' \rangle + \tfrac{1}{2} \on{tr} \big[ A A' \big] \geq G (x, q', A'), \ \text{ for all } (A', q') \in \mathbb{S}^d \times \bR^d \Big\}
		\] 
		is compact and convex. 
		We set 
		\begin{align*}
			\overline{\Theta}_F (x) := \Big\{ (A, q, r) \in \mathbb{S}^d \times \bR^d \times \bR \colon \ & \langle q, q' \rangle + \tfrac{1}{2} \on{tr} [ A A' ] + r \geq F (x, q', A'), \ r \leq \mathfrak{g} (x) \\& \text{ for all } (A', q') \in \mathbb{S}^d \times \bR^d \Big\}.
		\end{align*}
		{\em Step 1:}
		Let \((A, q, r) \in \mathbb{S}^d \times \mathbb{R}^d \times \bR\) be such that 
		\[
		\langle q, q' \rangle + \tfrac{1}{2} \on{tr} [ A A' ] + r \geq F (x, q', A')
		\]
		for all \((A', q') \in \mathbb{S}^d \times \bR^d\). Fix some \((A', q') \in \mathbb{S}^d \times \bR^d\). By \eqref{eq: F-G} and the positive homogeneity of~\(G\), for every \(t > 0\), 
		\[
		t \, \big( \langle q, q' \rangle + \tfrac{1}{2} \on{tr} [ A A' ] \big) + r \geq F (x, t q', t A') \geq t \, G (x, q', A') - \mathfrak{g} (x).
		\]
		Dividing by \(t\) and letting \(t \to \infty\) yields 
		\[
		\langle q, q' \rangle + \tfrac{1}{2} \on{tr} [ A A' ] \geq G (x, q', A'),
		\]
		which shows that \((A, q) \in \overline{\Theta} (x)\). 
		Taking \((A', q') = (0, 0)\) and using \eqref{eq: F-G} together with \(G (x, 0, 0) = 0\), we further get 
		\[
		r \geq F (x, 0, 0) \geq - \mathfrak{g} (x).
		\]
		Consequently,
		\begin{equation} \label{eq: alternative Theta F}
			\begin{split}
				\overline{\Theta}_F (x) = \Big\{ (A, q, r) \in \overline{\Theta} (x) \times [- \g (x), \g (x) ] \colon \ & \langle q, q' \rangle + \tfrac{1}{2} \on{tr} [ A A' ] + r \geq F (x, q', A')\\& \text{ for all } (A', q') \in \mathbb{S}^d \times \bR^d \Big\},
			\end{split}
		\end{equation}
		which in particular implies that \(\overline{\Theta}_F (x)\) is compact and convex. Conversely, \((A, q) \in \overline{\Theta} (x)\) implies \((A, q, \g (x)) \in \overline{\Theta}_F (x)\), as 
		\[
		\langle q, q' \rangle + \tfrac{1}{2} \on{tr} \big[ A A' \big] + \g (x) \geq G (x, q', A') + \g (x) \geq F (x, q', A'), 
		\] 
		by \eqref{eq: F-G}.
		
		{\em Step 2:}
		Next, similarly to the proof of \cite[Proposition~1]{Kusa_17}, we show that 
		\begin{align} \label{eq: rep F}
			F (x, q', A') = \inf_{(A'',\, q'',\, r) \in \overline{\Theta}_F (x)} \Big\{  \langle q', q'' \rangle + \tfrac{1}{2} \on{tr} [ A' A'' ] + r \Big\}.
		\end{align}
		The inequality \(\leq\) follows by the definition of \(\overline{\Theta}_F (x)\). For the converse direction, take \((A, q) \in \mathbb{S}^d \times \bR^d\). As \((A', q') \mapsto F (x, q', A')\) is a proper concave function with domain \(\mathbb{S}^d \times \bR^d\), its superdifferential at \((A, q)\) is non-empty (see \cite[Theorem~23.4, p.~217]{Rocka_CA}), i.e., there exists \((A'', q'') \in \mathbb{S}^d \times \bR^d\) such that, for all \((A', q') \in \mathbb{S}^d\times \bR^d\), 
		\begin{equation} \label{eq: supergradient}
			\begin{split} 
				F (x, q', A') &\leq F (x, q, A) + \langle q' - q, q'' \rangle + \tfrac{1}{2} \on{tr} \big[ (A' - A) A'' \big]
				\\&= \langle q', q'' \rangle + \tfrac{1}{2} \on{tr} \big[A' A'' \big] + r, 
			\end{split}
		\end{equation}  
		with \(r := F (x, q, A) - \langle q, q'' \rangle - \frac{1}{2} \on{tr} [ A A'' ]\). Step 1 shows that \((A'', q'') \in \overline{\Theta} (x)\). Furthermore, using this fact, we obtain that 
		\begin{align*}
			r \leq F (x, q, A) - G (x, q, A) \leq \g (x).
		\end{align*}
		Consequently, \((A'', q'', r) \in \overline{\Theta}_F (x)\) and, together with \eqref{eq: supergradient}, the inequality \(\geq\) in \eqref{eq: rep F} is proved. 
		
		\textit{Step 3:} 
		By \cite[Proposition~1]{Kusa_17}, for every \((A', q') \in \mathbb{S}^d \times \bR^d\), we have 
		\begin{align*} 
			\inf_{(A, q) \in \overline{\Theta} (x)} \Big\{ \langle q, q' \rangle + \tfrac{1}{2} \on{tr} \big[ AA' \big] \Big\} = G (x, q', A').
		\end{align*} 
		Together with \eqref{eq: rep F}, we obtain, for all \((A', q') \in \mathbb{S}^d \times \bR^d\) and \(s \in \bR\),
		\[
		\inf_{(A,q, r) \in \overline{\Theta}_F (x)} \Big\{  \langle q, q' \rangle + \tfrac{1}{2} \on{tr} [ A A' ] + sr \Big\}= 
		\begin{cases}
			s \, F (x, q' / s, A' / s), & s > 0, \\
			G (x, q', A'), & s = 0, \\
			G (x, q', A') + s \, \mathfrak{g} (x), & s < 0. 
		\end{cases}
		\]
		While the first case is immediate, the last two cases follow from the fact that, for all \((A, q) \in \overline{\Theta} (x)\), we have \((A, q, \g (x)) \in \overline{\Theta}_F(x)\), as shown at the end of Step 1. 
		All of these expressions are Borel measurable in \(x\). In view of \cite[Theorem~III.15]{CV_77}, this implies that \(x \mapsto \overline{\Theta}_F (x)\) is a measurable correspondence. Set \(\kappa_1 (x) := \sup_{\|A\| + \|q\| \leq 1} | G (x, q, A) |\) and \(\kappa_2 (x) := \mathfrak{g} (x)\), where we replace \(\kappa_i (x)\) by \(1\) in case it vanishes. Furthermore, define 
		\[
		\tilde{\Theta}_F (x) := \Big\{ \big( A / \kappa_1 (x), q / \kappa_1 (x), r / \kappa_2 (x) \big) \colon (A, q, r) \in \overline{\Theta}_F (x) \Big\}.
		\]
		As \(x \mapsto \overline{\Theta}_F (x)\) is Borel measurable, the same is true for \(x \mapsto \tilde{\Theta}_F (x)\) by virtue of \cite[Theorem~6.5]{himmelberg}. Furthermore, as \(\overline{\Theta}_F (x)\) has compact and convex values, the same is also true for \(\tilde{\Theta}_F (x)\). Moreover, it is straightforward to see that the values of \(\tilde{\Theta}_F\) are contained in a fixed ball whose radius only depends on \(d\). 
		
		{\em Step 4:}
		By the parametrization result \cite[Theorem~9.6.2]{AF_09}, there exists a compact Polish space \(\Gamma\) and measurable functions \((\tilde{a}, \tilde{b}, \tilde{g}) \colon \Gamma \times \bR^d \to \mathbb{S}^d \times \bR^d \times \bR\) such that
		\[
		\tilde{\Theta}_F (x) = \Big\{ (\tilde{a} (\gamma, x), \tilde{b} (\gamma, x), \tilde{g} (\gamma, x)) \colon \gamma \in \Gamma \Big\}, \quad x \in \bR^d,
		\] 
		and, for all \(x \in \bR^d, \gamma, \gamma' \in \Gamma\), 
		\[
		\| \tilde{a} (\gamma, x) - \tilde{a}(\gamma', x) \| + \| \tilde{b} ( \gamma, x) - \tilde{b} (\gamma', x) \| + | \tilde{g} (\gamma, x) - \tilde{g} (\gamma', x) | \leq c\, d_\Gamma (\gamma, \gamma'), 
		\] 
		where the constant \(c > 0\) only depends on the space dimension \(d\).
		Consequently, the triplet 
		\[
		\overline{a} (\gamma, x) := \tilde{a} (\gamma, x) \kappa_1 (x), \quad \overline{b} (\gamma, x) := \tilde{b} (\gamma, x) \kappa_1 (x), \quad \overline{g} (\gamma, x) := \tilde{g} (\gamma, x) \kappa_2 (x), 
		\] 
		satisfies \eqref{eq: for (a) 0}, \eqref{eq: properties g}, and 
		\[
		\overline{\Theta}_F (x) = \Big\{ (\overline{a} (\gamma, x), \overline{b} (\gamma, x), \overline{g} (\gamma, x)) \colon \gamma \in \Gamma \Big\}.
		\]
		Furthermore, assumption (iii) shows that, for every \(A \in \mathbb{S}^d_+, x \in \bR^d\) and \(\gamma \in \Gamma\), 
		\[
		\tfrac{1}{2} \on{tr} \big[ \overline{a} (\gamma, x) A \big] \geq \inf_{(A^*, p^*) \in \overline{\Theta} (x)} \tfrac{1}{2} \on{tr} \big[ A^* A \big] = G (x, 0, A) \geq 0, 
		\]
		which implies that \(\overline{a} (\gamma, x) \in \mathbb{S}^d_+\). 
		Summarizing, we conclude that \(F\) has the claimed HJB structure with coefficients \((\overline{a}, \overline{b}, \overline{g})\).
		The remaining claims concerning ellipticity and the \(L_d\)-drift bound only involve elements of \(\overline{\Theta} (x)\) and follow from straightforward computations. 
\qed 
 
 {\small
\bibliographystyle{siam}
\bibliography{references}}

\end{document}